\documentclass[final]{siamltex}
\usepackage{cite}
\usepackage{graphicx,bbm,pstricks}
\usepackage{subfigure}
\usepackage{amsmath,amsfonts,epsfig,amsbsy}
\newcommand{\RARR}[3]{#1
  \;\displaystyle\mathop{\displaystyle\longrightarrow}^{#3}\; #2}
\newcommand{\RARRlong}[3]{#1
  \;\displaystyle\mathop{-\!\!\!-\!\!\!-\!\!\!-\!\!\!-\!\!\!\!\displaystyle
  \longrightarrow}^{#3}\; #2}

\newcommand{\LRARR}[4]{{\mbox{ \raise 0.4 mm \hbox{$#1$}}} \;
  \mathop{\stackrel{\displaystyle\longrightarrow}\longleftarrow}^{#3}_{#4}
  \; {\mbox{\raise 0.4 mm\hbox{$#2$}}}}

\newcommand{\vecx}{{\mathbf x}}
\newcommand{\vecy}{{\mathbf y}}
\newcommand{\vecz}{{\mathbf z}}
\newcommand{\vecq}{{\mathbf q}}
\newcommand{\vecr}{{\mathbf r}}
\newcommand{\vecX}{{\mathbf X}}
\newcommand{\vecv}{{\mathbf v}}
\newcommand{\vecn}{{\mathbf n}}
\newcommand{\vecp}{{\mathbf p}}
\newcommand{\cT}{{\mathcal T}}
\newcommand{\dt}{{\mbox{d}t}}
\newcommand{\dx}{{\mbox{d} \vecx}}
\newcommand{\boldnu}{{\boldsymbol \nu}}
\newcommand{\er}{{\mathbb R}}

\newcommand{\divergence}{\mathop{\mbox{div}}}
\newcommand{\picturesAB}[6]{
\centerline{
\hskip #4
\raise #3 \hbox{\raise 0.9mm \hbox{(a)}}
\hskip #5
\epsfig{file=#1,height=#3}
\hskip #6
\raise #3 \hbox{\raise 0.9mm \hbox{(b)}}
\hskip #5
\epsfig{file=#2,height=#3}
}}
\newcommand{\picturesCD}[6]{
\centerline{
\hskip #4
\raise #3 \hbox{\raise 0.9mm \hbox{(c)}}
\hskip #5
\epsfig{file=#1,height=#3}
\hskip #6
\raise #3 \hbox{\raise 0.9mm \hbox{(d)}}
\hskip #5
\epsfig{file=#2,height=#3}
}}

\title{Adaptive Finite Element Method Assisted by
Stochastic Simulation of Chemical Systems}
\author{Simon L. Cotter\thanks{Mathematical Institute, University
of Oxford, 24-29 St. Giles', Oxford, OX1 3LB, United Kingdom
(cotter@maths.ox.ac.uk; erban@maths.ox.ac.uk).}
\and
Tom\'a\v s Vejchodsk\'y\thanks{Institute of Mathematics,
Czech Academy of Sciences, \v Zitn\' a 25, 115 67 Praha 1,
Czech Republic (vejchod@math.cas.cz).}
\and
Radek Erban$^*$}

\begin{document}
\maketitle
\begin{abstract}
Stochastic models of chemical systems are often analysed by solving
the corresponding Fokker-Planck equation which is a drift-diffusion
partial differential equation for the probability distribution function.
Efficient numerical solution of the Fokker-Planck equation
requires adaptive mesh refinements. In this paper, we present
a mesh refinement approach which makes use of a stochastic
simulation of the underlying chemical system.
By observing the stochastic trajectory for a relatively short
amount of time, the areas of the state space with non-negligible
probability density are identified. By refining the finite
element mesh in these areas, and coarsening elsewhere, a suitable
mesh is constructed and used for the computation of the probability
density.
\end{abstract}

\section{Introduction}
Stochastic simulation algorithms (SSAs) have been successfully used in
recent years to understand a number of biochemical models
\cite{Arkin:1998:SKA,Kar:2009:ERN,Villar:2002:MNR}. However, a
systematic analysis of these models is challenging because of the
computational intensity of SSAs. A suitable alternative to stochastic
simulation is a solution of the chemical Fokker-Planck equation
\cite{Gillespie:2000:CLE,Erban:2009:ASC}.  Consider a well-mixed
system of N chemical species and denote by ${\bf x} = (x_1, \ldots,
x_N)$ a vector of concentrations of these species. The stationary
chemical Fokker-Planck equation is a drift-diffusion partial
differential equation (PDE) for an $N$-dimensional probability
distribution function $p \equiv p(\vecx)$ where $\vecx \in \Omega
\subset \er^N$, which can be written in the following form:
\begin{equation}
\divergence
\big[
D(\vecx)
\nabla
p(\vecx)
-
\vecv(\vecx)
p(\vecx)
\big]
=
0
\label{SFPE}
\end{equation}
where
$D \equiv D(\vecx): \Omega \to \er^{N\times N}$ is the
diffusion matrix, $\vecv \equiv \vecv(\vecx): \Omega \to \er^N$
is the drift term and $\divergence$ is the divergence operator.
There have been several methods developed in the literature to solve
(\ref{SFPE}) for moderately large $N$. They include adaptive finite element
methods (FEMs) which are commonly used for $N \le 3$ \cite{Babuska:1978:EEA}
and sparse grid approaches which are applicable for larger values
of $N$ \cite{Zenger:1991:SG}. Adaptive FEMs can be used to identify
a suitable mesh which is refined in crucial regions and not in others.
Although these approaches can be used for the chemical
Fokker-Planck equation, they do not exploit the fact that its
solution is a probability distribution of a stochastic process.

In this paper, we present an adaptive mesh construction which is
suitable for the FEM solution of (\ref{SFPE}) if this equation arises
from modelling of stochastic chemical systems. The main idea is to
exploit the fact that stochastic trajectories spend a significant
amount of time in parts of the state space where the mesh refinement
is needed. Since adaptive FEMs are mostly applicable for systems up
to $N=3$, we will focus on systems of 3 chemical species. However, the
presented methodology can be modified for larger (multiscale) chemical
systems provided that they have up to $n \le 3$ important (slow)
variables. In \cite{Erban:2006:GRN,Cotter:2011:CAM}, a method for
estimation of coefficients of an effective Fokker-Planck equation is
presented. This effective equation is of the same form as (\ref{SFPE})
but it is written in the dimension $n$ which is smaller than a total
number $N$ of chemical species. If one has a suitable SSA for
simulating the low dimensional slow dynamics
\cite{Cao:2005:SSS,Cao:2005:MSS,E:2005:NSS}, the presented mesh
refinement can be applied. However, since the main aim of this paper
is to present this novel numerical methodology, we will not consider
any dimensional reduction and restrict to systems which are directly
written in terms of $N=3$ chemical species.  In the following paper
\cite{Cotter:2011:SB}, we will show how this mesh refinement can be
used to study bifurcation behaviour of stochastic chemical systems.

The paper is organised as follows. In Section \ref{sec:CFPE}, we
introduce our notation, the Gillespie SSA and the chemical
Fokker-Planck equation. In Section \ref{sec:FEM}, we introduce the
standard finite element framework that we will use throughout the
paper. We also formulate the problem in its weak form in order to
define the problem to be solved. In Section \ref{sec:saFEM}, we
introduce the stochastic simulation assisted adaptive Finite Element
Method (saFEM). In Section \ref{sec:param}, we offer some insight into
how the algorithmic parameters of the saFEM may be chosen in
practise. In Section \ref{sec:Imp}, implementational issues related to
the algorithm are addressed. Numerical results concerning convergence
of the method are presented in Section \ref{sec:res}.

\section{Chemical Fokker-Planck Equation}
\label{sec:CFPE}

Let us consider a well-mixed system of $N$ chemical species $X_i$, $i
= 1, 2, \dots, N$, which are subject to $M$ chemical reactions $R_k$,
$k=1, 2, \dots, M$.  The state of the system is described by the state
vector $\vecX(t) = [X_1(t), X_2(t), \dots, X_N(t)]$ where $X_i(t)$
denotes the number of molecules of the corresponding chemical species.
Each chemical reaction is described by its propensity function and the
stoichiometric vector \cite{Gillespie:1977:ESS,Erban:2007:PGS}.  The
propensity function $\alpha_k(\vecx)$ is defined in such a way that
$\alpha_k(\vecx) \dt$ is the probability that the $k$-th reaction
occurs in the infinitesimally small interval $[t,t+\dt)$ provided that
$\vecX(t) = \vecx.$ The stoichiometric vector is $\boldnu_k =
[\nu_{k1}, \nu_{k2}, \dots, \nu_{kN}]$ where $\nu_{ki}$ is the change
in $X_i$ during reaction $R_k$. Stoichiometric vectors form the
corresponding stoichiometric matrix $\nu = (\nu_{ki})_{k,i=1}^{M,N}$.

The time evolution of the state vector $\vecX(t)$ is often simulated
by the Gillespie SSA \cite{Gillespie:1977:ESS} which is described
in Table \ref{SSAtable}.
\begin{table}
\framebox{%
\hsize=0.97\hsize
\vbox{
\leftskip 6mm
\parindent -6mm
{\bf [1]} Calculate propensity functions $\alpha_k(\vecX(t))$,
$k=1,2,\dots,M$.

\smallskip

{\bf [2]} Waiting time $\tau$ till next reaction is given by \eqref{eq:tau}.

\smallskip

{\bf [3]} Choose one $j\in\{1,2,\ldots,M\}$, with probability $\alpha_j(\vecX(t))/\alpha_0(\vecX(t))$,
and perform reaction $R_j$, by adding $\nu_{ji}$ to each
$X_i(t)$ for all $i = 1,2,\dots, N$.

\smallskip

{\bf [4]} Continue with step [{\bf 1}] with time $t=t+\tau$.

\par \vskip 0.8mm}
}\vskip 1mm
\caption{{\it The pseudo code for the Gillespie SSA.}
\label{SSAtable}}
\end{table}
Given the values of the propensity functions (step [{\bf 1}]),
the waiting time to the next reaction is given by:
\begin{equation}\label{eq:tau}
  \tau = -\frac{\log\left(u\right )}{\alpha_0(\vecX(t))},
  \quad \mbox{where} \quad \alpha_0(\vecX(t)) = \sum_{k=1}^M \alpha_k(\vecX(t)),
\end{equation}
and $u$ is a uniformly distributed random number in $(0,1)$.
The reaction $R_j$ is chosen in step [{\bf 3}] using another
uniformly distributed random number.

The stationary probability distribution corresponding to the
chemical system can be approximated by solving the chemical
Fokker-Planck equation \cite{Gillespie:2000:CLE}. This equation
can be written in the form (\ref{SFPE}) where the
diffusion and drift coefficients are:
\begin{eqnarray}\label{eq:CFPE}
d_{ij}(\vecx) &=& \frac{1}{2}
     \sum_{k=1}^M \nu_{ki} \nu_{kj} \alpha_k(\vecx), \qquad
     i,j = 1, 2, \dots, N, \\ \label{eq:CFPE2}
v_{i}(\vecx) &=& \sum_{k=1}^M\nu_{ki}\alpha_k(\vecx)
- \sum_{j=1}^N \frac{\partial d_{ij}}{\partial x_j}(\vecx),
\qquad i = 1, 2, \dots, N.
\end{eqnarray}
This equation is solved on a bounded domain $\Omega \subset [0,\infty)^N$ in which the
vast majority of the invariant probability density sits \cite{Erban:2009:ASC}.
On the boundary $\partial\Omega$ we introduce the homogeneous
Neumann boundary condition:
\begin{equation}
\label{eq:strong2}
\bigl[ D(\vecx)\nabla p(\vecx) + \vecv(\vecx) p(\vecx) \bigr]
\cdot \vecn(\vecx) = 0,  \qquad\qquad \rm{for}\ \vecx \in \partial \Omega,
\end{equation}
where $\vecn(\vecx)$ is the outward facing normal at $\vecx
\in \partial \Omega$.  We seek a solution of (\ref{SFPE})
which corresponds to the probability distribution function.
Therefore, we impose the following normalisation condition:
\begin{equation}
\label{eq:strong1}
\int_\Omega p(\vecx) \,\mbox{d}\vecx = 1.
\end{equation}
The choice of the computational domain $\Omega$ might be problematic
if we have no a priori information about the problem (\ref{SFPE}) with
(\ref{eq:strong2}) and (\ref{eq:strong1}). In this case, the stochastic
simulations provide a reliable tool to determine the correct
$\Omega$ as we will show in Section~\ref{sec:2}. However, before
discussing the details of saFEM, we have to introduce some finite
element terminology. This will be done in the next section.

\section{Finite Element Method}\label{sec:FEM}

Equation (\ref{SFPE}) with boundary condition (\ref{eq:strong2}) can
be numerically solved by the FEM. Since the finite element formulation
is based on the corresponding weak formulation, we first introduce the
weak solution $p\in H^1(\Omega)$ by the equality
\begin{equation}
\label{eq:weak1}
a(p,\phi) = 0, \qquad\qquad \forall \phi\in H^1(\Omega).
\end{equation}
The  bilinear form $a(\cdot,\cdot)$ is naturally given by
\begin{equation}
\label{eq:blf}
a(p,\phi) = \int_\Omega (D(\vecx) \nabla p(\vecx) + p(\vecx) \vecv(\vecx))
\cdot \nabla\phi(\vecx) \, \mbox{d}\vecx.
\end{equation}

The finite element formulation is obtained by projecting the weak formulation
(\ref{eq:weak1}) into a finite dimensional subspace $V_h$ of $H^1(\Omega)$.
Thus, we seek $p_h \in V_h$ such that
\begin{equation}
\label{eq:FEM}
a(p_h,\phi_h) = 0, \qquad\qquad \forall \phi_h\in V_h.
\end{equation}
We note that taking a basis $\phi_1,\phi_2,\dots,\phi_m$ of $V_h$,
we can express the finite element solution as
$$
  p_h(\vecx) = \sum\limits_{i=1}^m p^i \phi_i(\vecx),
  \quad \vecx\in\Omega.
$$
Here, the coefficients $p^i \in \er$ solve the system of
linear algebraic equations
\begin{equation}
\label{eq:linsys}
A\vecp = 0,
\end{equation}
where $\vecp = (p^1,p^2,\dots,p^m)^T$
and the stiffness matrix $A\in \er^{m\times m}$
is defined by its entries
\begin{equation}
\label{eq:stiff}
A_{ij} = a(\phi_j,\phi_i), \quad i,j=1,2,\dots,m.
\end{equation}

We construct the finite element space $V_h$ and the corresponding
basis functions in the standard way \cite{Ciarlet:1978:FEM}. In
what follows, we will focus on computations in three dimensions,
i.e. $N=3$. We consider the finite element mesh $\cT_h$ consisting
of tetrahedral elements. The lowest-order finite element space $V_h$
then consists of globally continuous and piecewise linear functions
over the mesh $\cT_h$:
\begin{equation}
\label{eq:Vh}
  V_h = \{ \phi_h \in H^1(\Omega) : \phi_h|_K \in \mathbb{P}^1(K)
  \mbox{ for all } K\in\cT_h \},
\end{equation}
where $\mathbb{P}^1(K)$ stands for the space of linear functions over the
tetrahedron $K \in \cT_h$.
If $\vecq_j \in \er^3$, $j=1,2,\dots,m$, stand for the nodes of
the mesh $\cT_h$
then the standard finite element basis functions $\phi_i$ are uniquely
determined by the condition
$$
  \phi_i(\vecq_j) = \delta_{ij}, \quad i,j=1,2,\dots,m,
$$
where $\delta_{ij}$ stands for Kronecker's symbol.

An efficient solution to problem (\ref{eq:FEM}) can be obtained by
employing adaptively refined meshes \cite{Verfurth:1996:RAP}.
The optimally adapted mesh leads
to an approximation with the smallest error, provided the number of
degrees of freedom is fixed.  Practically, the optimal mesh can be
hard to determine, but meshes close to the optimal can be found.
These meshes are fine in regions where the solution exhibits steep
gradients, boundary layers, interior layers or singularities, and
they are relatively coarse in the other regions.

The standard numerical approach for construction of nearly optimal
meshes is the adaptive algorithm based on suitable a posteriori error
estimators \cite{Babuska:1978:EEA,Verfurth:1996:RAP}.  This
algorithm starts with an initial coarse mesh and refines it adaptively
by a sequence of refinement steps. In each step, problem
(\ref{eq:FEM}) has to be solved on the actual mesh, an error indicator
has to be computed for each element, and based on these indicators the
mesh is refined at suitable places.  Using the mesh refinements
assisted by stochastic simulations, we can avoid the sequence of
refinement steps and construct a suitably adapted mesh at once.

\section{Adaptive mesh refinement assisted by stochastic simulations}
\label{sec:saFEM}
Since the chemical Fokker-Planck equation \eqref{SFPE} is a continuous
approximation of the evolution of the probability density given by the
Gillespie SSA\cite{Gillespie:1977:ESS}, one can expect that
trajectories simulated from the SSA will be informative. Once the
trajectory has reached probabilistic equilibrium\footnote{Equilibrium
  can be reached simply by running the SSA until the state of the
  trajectory is sufficiently decorrelated from its starting position,
  or by starting at a steady state of the mean field approximation
  (\ref{eq:ode}) of the chemical system.}, regions surrounding the
path of the trajectory are likely to be regions with non-negligible
invariant density with respect to the steady state Fokker-Planck
equation. We should be aiming to refine the finite element mesh in
regions where the rate of change of the gradient of the invariant
density is larger. The area in which we should be trying to refine our
mesh can be well approximated by the region which has non-negligible
invariant density. Therefore, stochastic simulations of the chemical
system can be informative about a good choice of finite element mesh.

The construction of the locally adapted mesh assisted by stochastic
simulations is done in three stages.  In Stage I, we use the
stochastic simulations to identify those regions of the state space,
where the system spends most of its time, and hence where we are
interested to resolve the problem with the highest accuracy.  In Stage
II, we identify the computational domain $\Omega$ as a cuboid that
covers the region from Stage I and its neighbourhood.  In Stage III,
we construct the actual mesh in the computational domain
$\Omega$ based on the information from Stage I.  In what
follows, we provide detailed description of these three stages in the
case of $N=3$ chemical species. Following the notation of Section
\ref{sec:CFPE}, the chemical species will be denoted as $X_1$, $X_2$
and $X_3$. The mesh generation algorithm is summarised in Table
\ref{Meshgen Table}.

\begin{table}
\framebox{%
\hsize=0.97\hsize
\vbox{
\leftskip 6mm
\parindent -6mm
{\bf [1]} Identify steady states (stable and unstable)
of the mean field approximation of the system, or in
the case of oscillatory systems, one (or more) coordinates
along the limit cycle. We denote these $S\in \mathbb{N}$
points by $\{\vecy_k\}_{k=1}^S$.

\smallskip

{\bf [2]} Run one (or more) Gillespie simulations of length
$B+T>0$ for each of the initial conditions $\{\vecy_k\}_{k=1}^S$.
Here, $B \ge 0$ is the length of a possible initial transient
and $T>0$.  Take a subsample of $Q>0$ points per unit time
per trajectory in the time interval $[B,B+T]$ and denote them
as in (\ref{calAdef}). We also denote by $x_i^{\rm max}$
(resp. $x_i^{\rm min}$) the maximal (resp. minimal) value
of the $i$-th chemical species, $i = 1, 2, 3$, during time
intervals $[B,B+T]$ of all $S$ simulations.

\smallskip

{\bf [3]} Use (\ref{definitionofgamma}), to
define a neighbourhood $\Gamma \subset [0,\infty)^N$
of the points (\ref{calAdef}) as the
region of the domain in which we require the finest level of refinement
of the mesh. This is made up of ellipsoids around each point with
radii given by \eqref{eq:uncert} with parameter $\beta_1 > 0.$

\smallskip

{\bf [4]} Define the domain of solution $\Omega$ by
\eqref{eq:domain}, using parameter $\beta_2 > 0.$

\smallskip

{\bf [5]} Start with the mesh as one single cuboid covering
the whole of $\Omega$.

\smallskip

{\bf [6]} Loop over all cuboids in the mesh.
If the cuboid is sufficiently close (as given by \eqref{eq:cond})
to $\Gamma$, then split the cuboid into eight equally sized cuboids.
Update the list of hanging nodes/hang type (face/edge hanging node
as shown in Figure \ref{fi:hnodes} (left)).

\smallskip

{\bf [7]} Repeat step {\bf [6]} until the maximum resolution has
been reached after $H\in\mathbb{N}$ iterations, or the maximum
total number of cuboids has been reached.
\par \vskip 0.8mm}
}
\caption{{\it Mesh generation using SSA trajectories.}\label{Meshgen Table}}
\end{table}

\subsection{Stage I: Stochastic simulation}\label{sec:1}

In step {\bf [1]} of Table \ref{Meshgen Table}, we identify structures
that a priori should exist in the probability density. An indication
of the regions which may contain significant amounts of invariant density
can be found by analysing the 3-dimensional system of
ordinary differential equations (ODEs)
\begin{equation}\label{eq:ode}
\frac{\mbox{d}x_i}{\dt} = v_i(x_1,x_2,x_3), \qquad i \in \{1,2,3\},
\end{equation}
which is closely related the mean field approximation of the chemical
system. Steady states of ODE system (\ref{eq:ode}) will often coincide
with regions of the solution of the steady state Fokker-Planck
equation \eqref{SFPE} that have large density.  The most important
things to identify are stable steady states of \eqref{eq:ode}, which
can be found by solving an algebraic system $[v_1(x_1,x_2,x_3),v_2(x_1,x_2,x_3),v_3(x_1,x_2,x_3)] = 0$. In
oscillating systems, one can identify limit cycles. There are several
tools in the literature for analysis of ODEs of the form
(\ref{eq:ode}), such as AUTO \cite{Doedel:1997:ACB}.

In step {\bf [2]}, these steady states can then be used as starting
points for SSA trajectories, once the numbers of molecules of each
species have been rounded to the nearest integer. One might
additionally wish to ensure that each chain is starting in
probabilistic equilibrium by running a short simulation of length $B
\ge 0$ from these points before starting the sampling procedure. The
Gillespie SSA is detailed in Table \ref{SSAtable}. In the case of
limit cycles, several points along the limit cycle can be used as
starting points for the SSA. Either way, we denote the $S$ starting
positions for the trajectories by $\{\vecy_k\}_{k=1}^{S} \subset
\er^3$.  The trajectories simulated up to some time $B+T>0$ using the
SSA can help inform us about the regions of domain which will have
non-negligible invariant density.  Since the trajectories will contain
many points, we cannot store all of the points that are
simulated. Instead, we subsample from the trajectories at equidistant
time points, at a rate $Q > 0$ points per unit time. This leaves us
with a set of sampled points
\begin{equation}
\left\{\vecz_{k,l}\right\}_{k=1, l=1}^{S,\lfloor QT \rfloor} \subset \er^3
\label{calAdef}
\end{equation}
from the invariant distribution, where $\lfloor QT \rfloor$ denotes
the integer part of the real number $QT$. We hope to recover, from
this set of points (\ref{calAdef}), information about what might be an
optimal finite element mesh. We also denote by $x_i^{\rm max}$
(resp. $x_i^{\rm min}$) the maximal (resp. minimal) value of the
$i$-th chemical species, $i = 1, 2, 3$, during time intervals
$[B,B+T]$ of all $S$ simulations. These numbers will be useful in
(\ref{xranges})--(\ref{eq:uncert}).  In step {\bf [3]}, we make the
approximation that the region which contains the majority of the
invariant density, is a subset of the union of a set of ellipsoids
centred at each point (\ref{calAdef}), namely
\begin{equation}
\Gamma \equiv
\bigcup_{k,l} \mathcal{E}_{\vecr}(\vecz_{k,l}).
\label{definitionofgamma}
\end{equation}
Here, $\mathcal{E}_{\vecr}(\vecz_{k,l})$ is an ellipsoid with radii
${\vecr} = [r_1, r_2, r_3]$ centred at point $\vecz_{k,l}$ for $k=1,
2, \dots, S,$ $l=1, 2, \dots, \lfloor QT \rfloor$. These radii can be
picked to be proportional to the range of each chemical species using
the parameter $\beta_1 > 0$ and numbers $x_i^{\rm min}$,
$x_i^{\rm max}$, $i = 1, 2, 3$, computed in step {\bf [2]}. Namely, we
define
\begin{equation}
\label{xranges}
x_i^{\rm range} = x_i^{\rm max}-x_i^{\rm min},
\qquad\qquad
i = 1, 2, 3,
\end{equation}
and
\begin{equation}
\label{eq:uncert}
{\bf r} = \beta_1 \left (x_1^{\rm range}, \,
x_2^{\rm range}, \, x_3^{\rm range} \right).
\end{equation}
One should note that the parameters B and T might in principle be
different for different initial conditions $y_k$, but to simplify the
notation, we do not stress this fact by using the notation $B_k$ and
$T_k$ and use simply $B$ and $T$.
\subsection{Stage II: Automatic detection of the domain}\label{sec:2}

Next, we would like to identify the domain $\Omega$ on which we wish
to solve \eqref{SFPE}. Since we already have an approximation of the
region which contains the majority of the probability density, we can
simply fit a cuboid around this region, and then extend it by a
factor. The only other condition that we enforce is that the domain
must be contained by the positive quadrant of the state space. In step
{\bf [4]}, we pick a parameter $\beta_2 > 0$ to define how much we
would like to extend the cuboid:
\begin{equation}
\label{eq:domain}
\Omega =  {\mathcal A}_1 \times {\mathcal A}_2 \times {\mathcal A}_3,
\end{equation}
where
\begin{equation*}
{\mathcal A}_i
=
\left( \max \{0, x_i^{\rm min}- \beta_2 \, x_i^{\rm range}\},
x_i^{\rm max} + \beta_2 \, x_i^{\rm range}
\right),
\qquad \qquad i = 1, 2, 3.
\end{equation*}
We then wish to solve \eqref{SFPE} on the domain $\Omega$ with
boundary conditions (\ref{eq:strong2}) on $\partial \Omega$.

\subsection{Stage III: Construction of the adaptive mesh}

\begin{figure}[thb]
\centering
\includegraphics[width=0.45\textwidth]{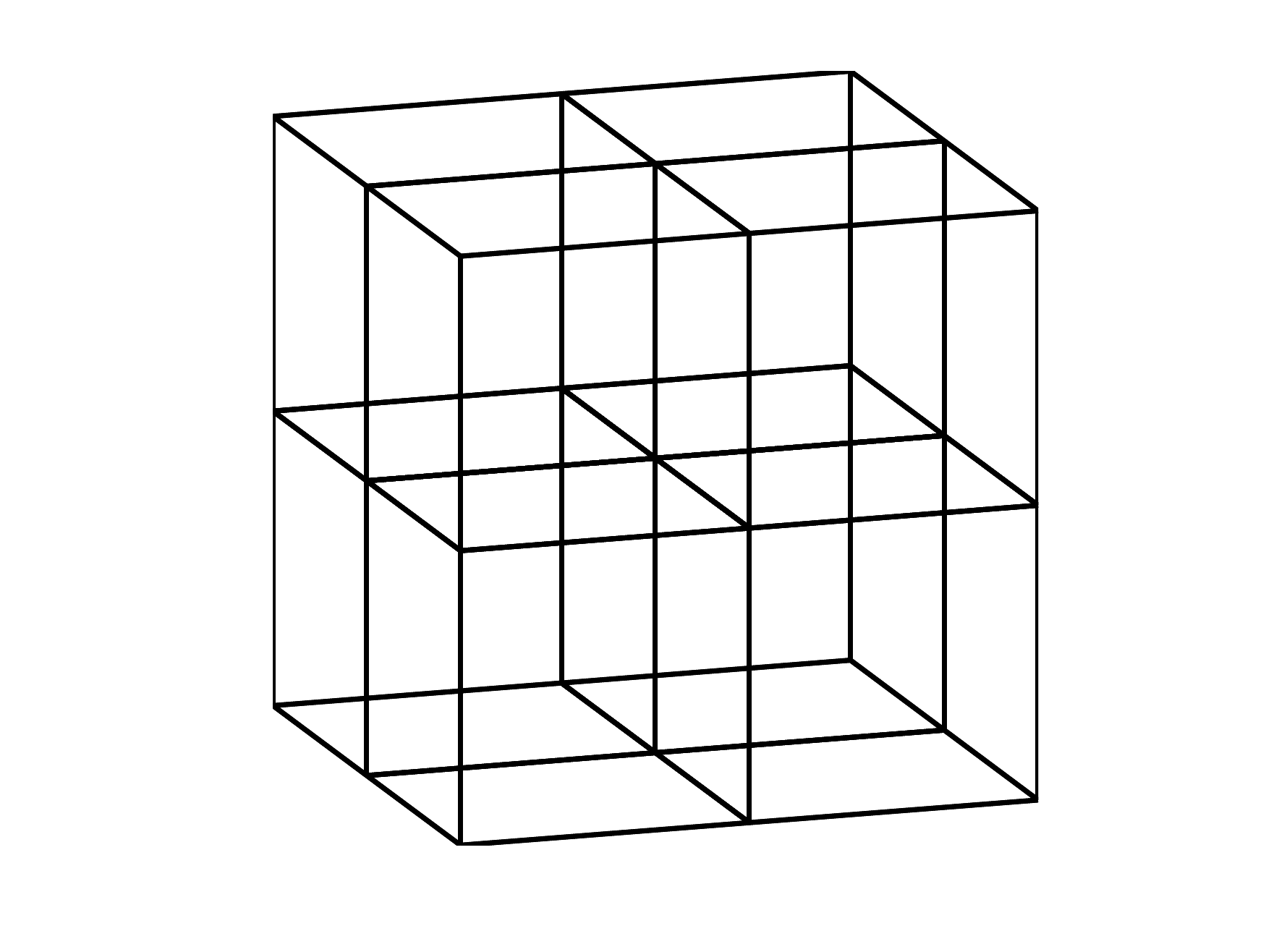}
\caption{{\it Example of a single cuboid after the proposed
refinement into $8$ sub-cuboids.}\label{split_cube}}
\end{figure}

In this stage, we construct the adaptively refined mesh.  In the
previous stage we naturally defined the computational domain as a
cuboid (\ref{eq:domain}).  Therefore, we base our refinement approach
on simple splitting of a cuboid into 8 congruent sub-cuboids, as shown
in Figure~\ref{split_cube}.  However, any other standard mesh
refinement technique can be utilised instead.

\begin{figure}[tb]
\picturesAB{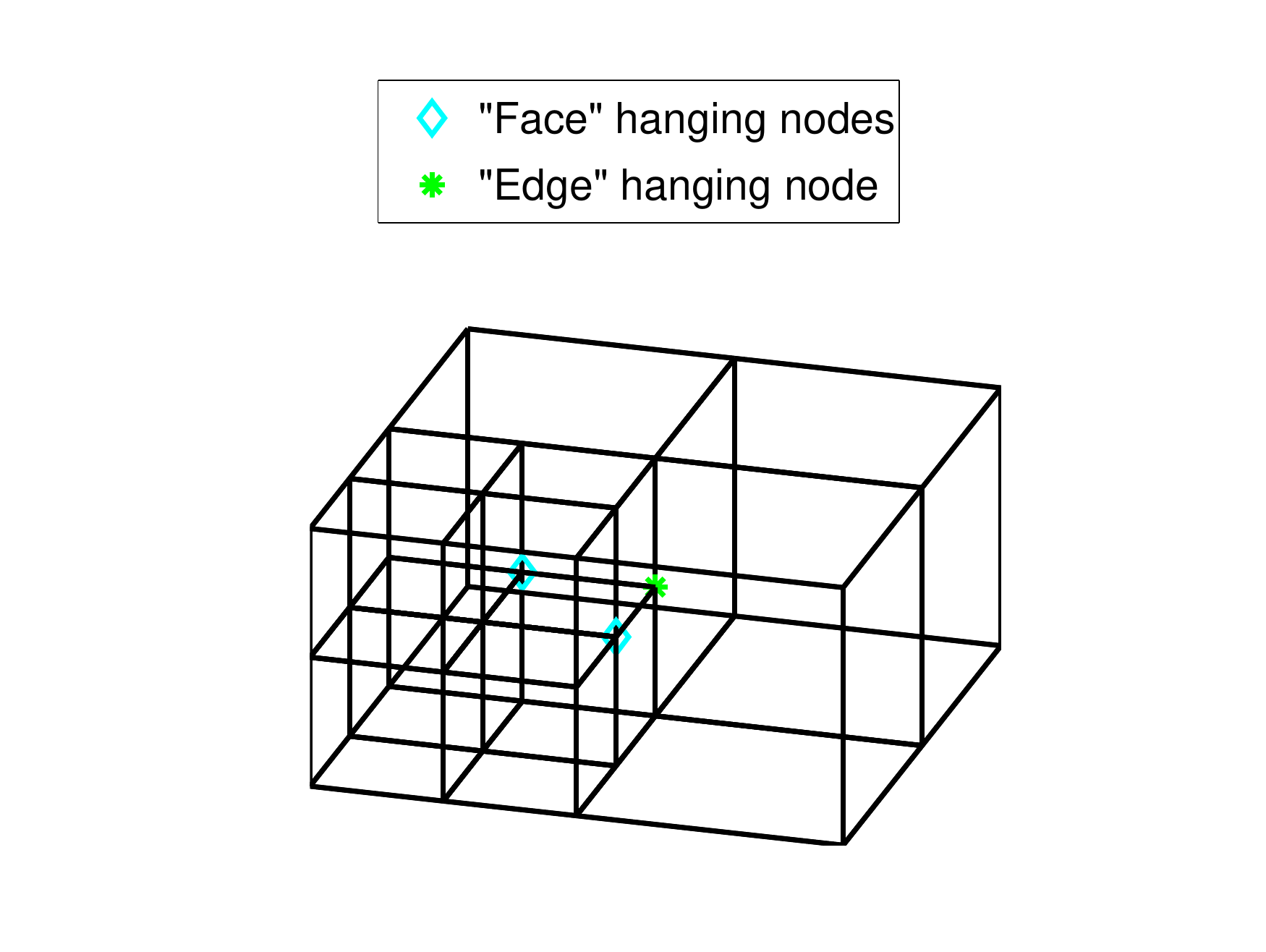}{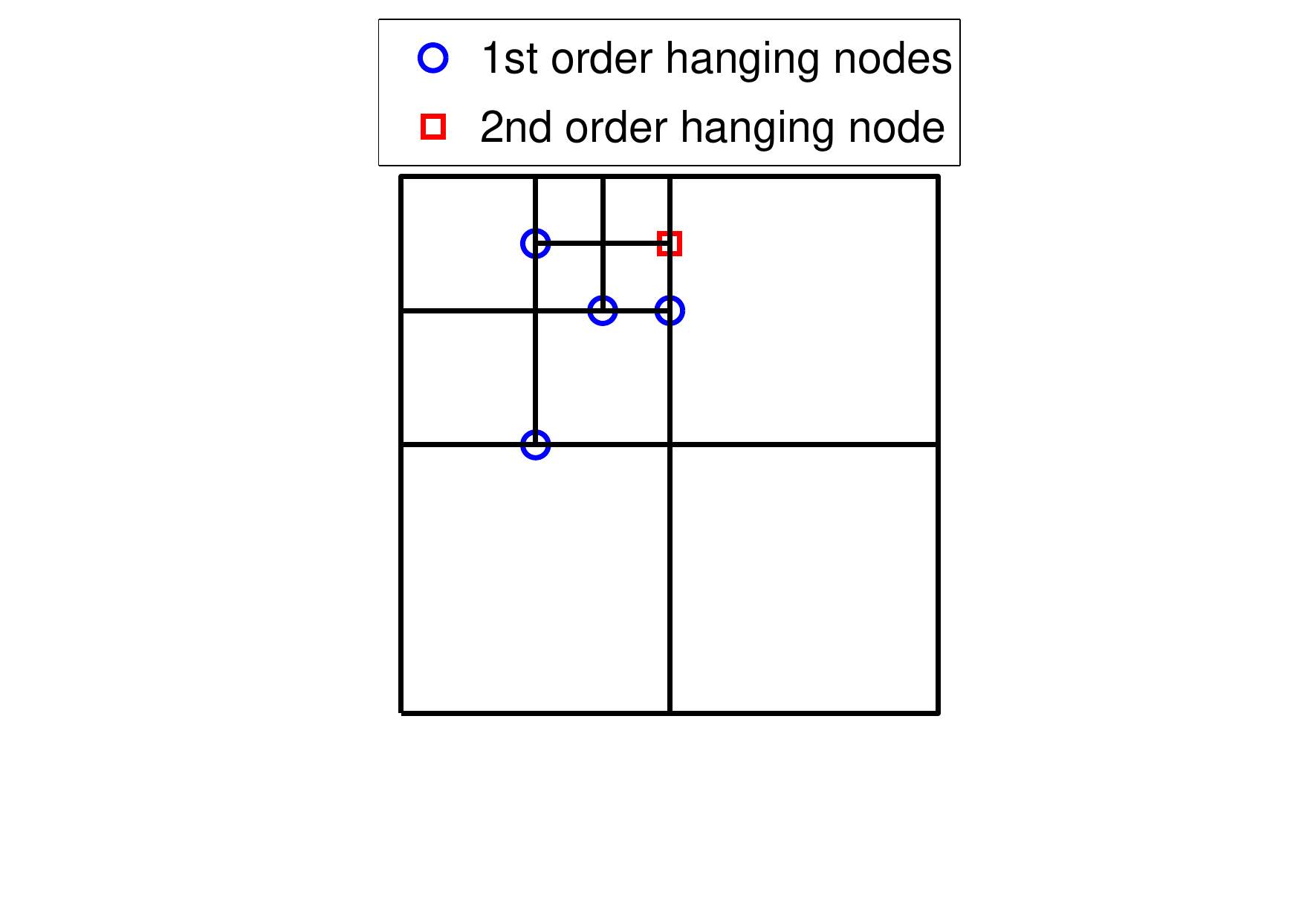}
{2.05in}{15mm}{-15mm}{10mm}
\caption{{\rm (a)} {\it $1^{st}$ order hanging nodes of different types in
cuboid meshes.} {\rm (b)} {\it
Hanging nodes of order $1$ and order $2$ on a {\rm 2D} mesh.}
\label{fi:hnodes}}
\end{figure}

An advantage of the refinement of cuboids into 8 sub-cuboids is its simplicity.
A disadvantage is that these refinements produce so-called hanging nodes
in the mesh, see Figure~\ref{fi:hnodes} (left).

In step {\bf [5]} of the mesh generation, our mesh consists of a
single large cuboid. We then, in step {\bf [6]}, iteratively refine
cuboids in the following fashion. In each iteration of the refinement
procedure, we loop over all of the current cuboids that exist
following the previous iteration of the method. For each of the
cuboids $I_1 \times I_2 \times I_3$, where $I_1,I_2,I_3 \subset
[0,\infty)$ are intervals, we compute the minimum distance between all
points in this cuboid, and the set of points (\ref{definitionofgamma})
in each coordinate, namely
\begin{equation*}
  \mbox{dist}(\Gamma,I_i)
  = \inf\{|a_i-b_i| : [a_1,a_2,a_3]\in \Gamma, \, b_i \in I_i
  \}, \qquad \qquad i = 1, 2, 3.
\end{equation*}
We refine the given cuboid $I_1 \times I_2 \times I_3$ if
\begin{eqnarray}
\label{eq:cond}
  \mbox{dist}(\Gamma,I_i) < |I_i|
  \qquad  \mbox{is satisfied for every} \; i = 1, 2, 3.
\end{eqnarray}
If a cuboid is to be refined, it is split into 8 cuboids of equal
volume. As detailed in step {\bf [7]}, this refinement condition can
be iterated for a set number of times $H\in\mathbb{N}$, where
$\mathbb{N} \equiv \{1,2,3,\ldots \}$, or until the number of cuboids
present in the mesh is as large as we would like or can deal with given
our computational resources.

Notice that the algorithm described above produces hanging nodes of
order one only. This means that two neighbouring elements in the mesh
are either of equal size or one of them has eight times greater volume
than the other one.  This is important for practical implementation,
because the hanging nodes of order one are much easier to work with
than the hanging nodes of higher orders. See \cite{Solin:2008:ALH} for
more details and Figure~\ref{fi:hnodes}(b) for an illustration in
2D. In order to guarantee the continuity of the approximation, the
value of the function at the hanging node is necessarily decided by
the values of the function at the vertices of the less refined
cube. Therefore the hanging nodes are not in actual fact additional
degrees of freedom in the problem.  Thus, the dimension of problem
(\ref{eq:linsys}) is equal to the total number of vertices in the
generated mesh, less the number of hanging nodes.  This special
treatment of hanging nodes actually enforces on us a mesh which is
highly refined in the regions which we wish it to be, and then the
mesh becomes gradually coarser as we move away from those regions.

\subsection{Tetrahedral mesh}

Once the cuboid mesh has been generated by steps {\bf [1]--[7]} of
Table \ref{Meshgen Table}, we can then
implement a finite element method on it. In the numerics shown in this
paper, we further split each cuboid into 6 path tetrahedra. However
there is nothing to say that we could not use cubic elements, but we
use tetrahedra to simplify some implementational issues.
\begin{figure}[tb]
\picturesAB{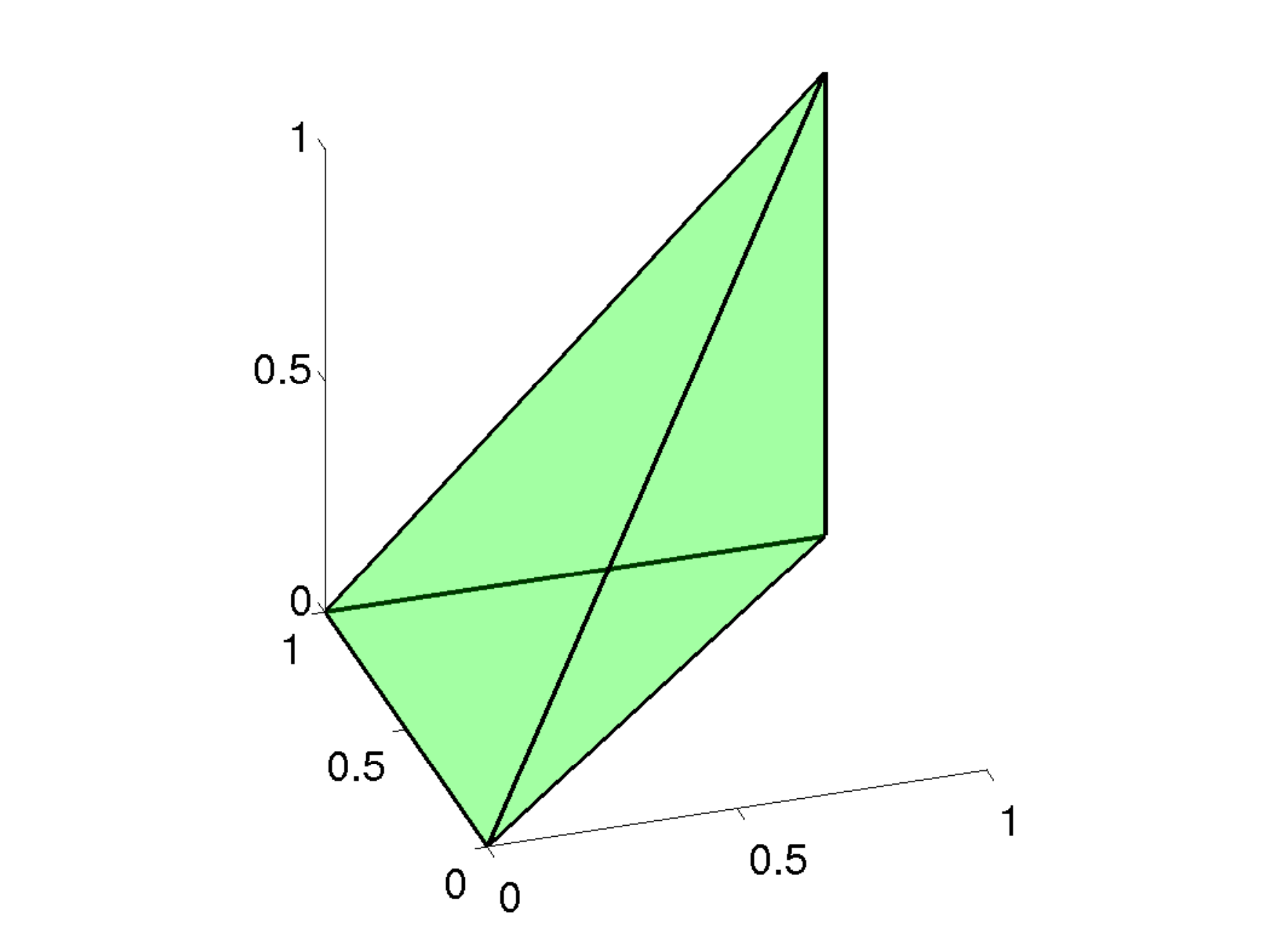}{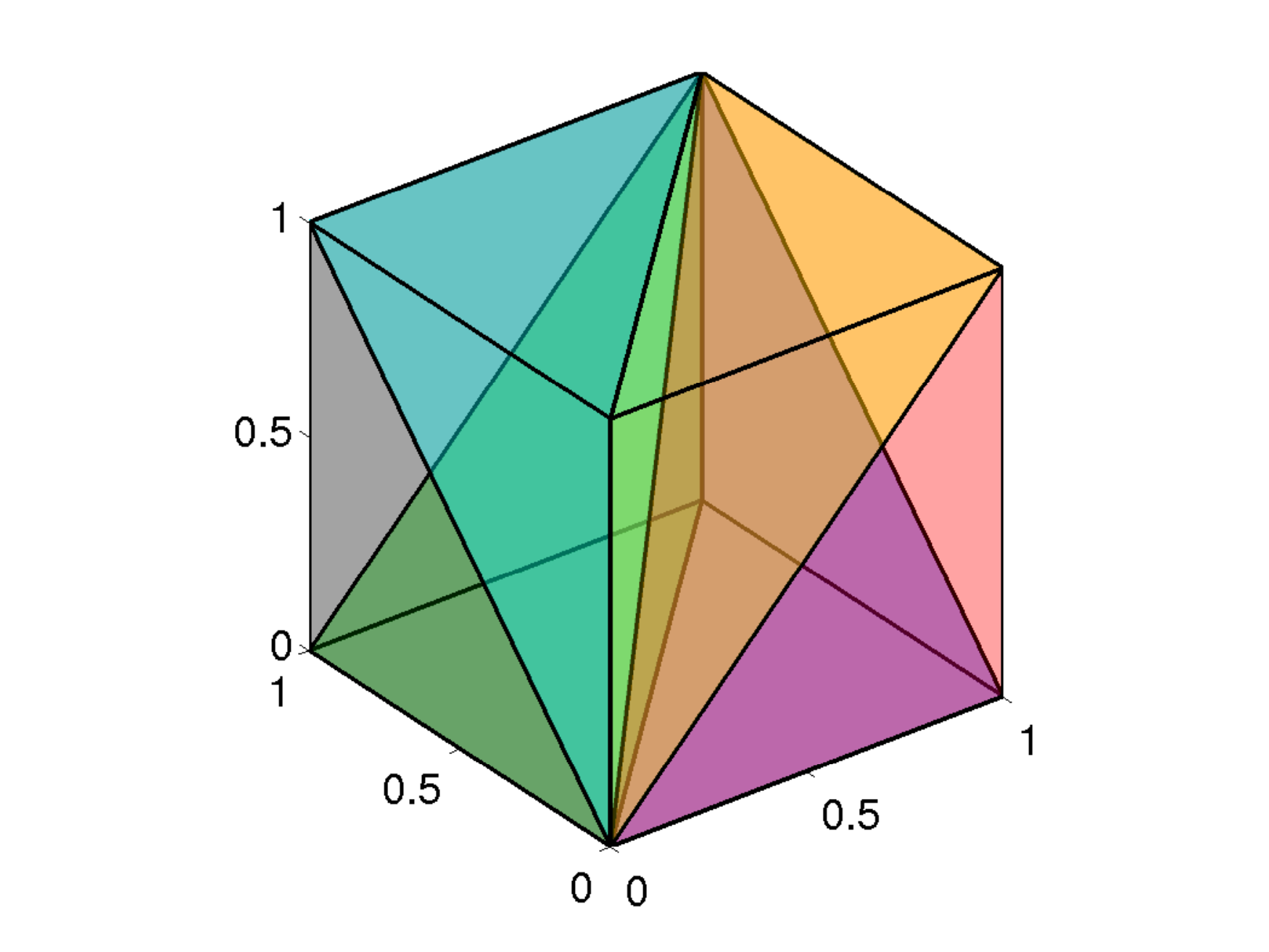}
{2.05in}{15mm}{-15mm}{10mm}
\caption{{\rm (a)} {\it Example of one element from the refinement scheme.}
{\rm (b)} {\it Example of a unit cube split into 6 elements of equal
  size and shape.}\label{elements}}
\end{figure}

We choose a refinement regime in which the elements are clustered in
groups of 6 non-obtuse\footnote{All six dihedral angles between its
  faces are less than or equal to $\pi/2$} tetrahedra which together
form each cuboid \cite{Beilina:2005:NTP}. Figure
\ref{elements}(a) shows an example element from the refinement
scheme. The idea is that if we map each cuboid to the unit cube, then
each element is exactly of the same shape and size as all the
others. Figure \ref{elements}(b) shows how 6
of these elements can tessellate into a cube.

A lot of refinement methods in the literature involve splitting each
element into several smaller elements\cite{Korotoc:2002:ATR}. However,
if this is not done in a clever way, this can lead to degradation of
the quality of the elements. That is, some of the angles of the
elements may become too small, leading to long thin elements, which
can lead to less accuracy in the approximation\cite{Bey:1995:TGR}.

\begin{figure}[thb]
\centering
\includegraphics[width=0.45\textwidth]{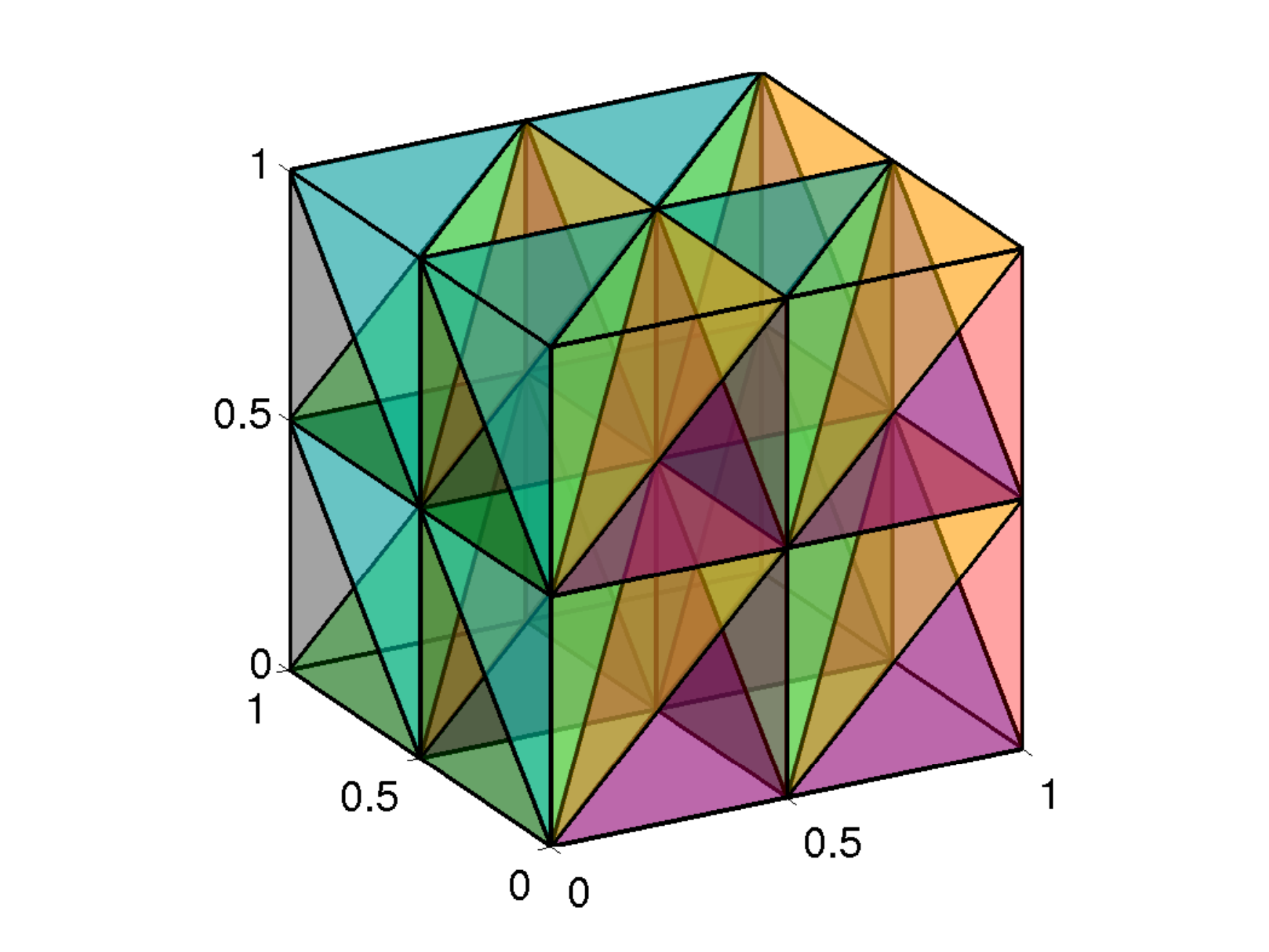}
\caption{\it Example of a unit cube after a refinement iteration.\label{cube_refined}}
\end{figure}

In the method that we propose here, instead of splitting each element,
we split each cuboid into 8 equally sized cuboids. Each cuboid is then
split, as before, into 6 equally sized and shaped (after a linear
transformation if not a cube) tetrahedral elements. Figure
\ref{cube_refined} shows a cube which has been refined once and been
split into 48 elements, with 27 vertices.

\section{Parameter Selection}\label{sec:param}
The algorithm has several parameters $S$, $T$, $Q$, $B$, $\beta_1$,
$\beta_2$ and $H$ whose values must be decided by the user. The
parameters $S$, $T$, $Q$, $B$ relate to the implementation of the
Gillespie SSA, and the remaining three $\beta_1$, $\beta_2$, $H$
relate to the selection of the domain and the mesh, given the recorded
output of the SSA. In this section we suggest how one might go about
selecting values for these parameters.

\subsection{SSA parameters}
The purpose of simulating the SSA trajectories, as described in
Subsection \ref{sec:1}, is to get an indication of the regions of the
domain which contain the majority of the invariant probability
density. There are four parameters to consider, $S \in {\mathbb N}$,
$T >0$, $Q > 0$ and $B \ge 0$. The parameter $B$ represents the length
of the simulated trajectory that we ignore at the beginning of the
simulation. This period of simulation is simply used to ensure that
the Markov chain has entered probabilistic equilibrium. This ensures
that it is highly likely the point in state space where the rest of
the trajectory starts has non-negligible density with respect to the
invariant distribution. The final parameter, $Q$, represents the rate
(per unit time) at which we record samples from each trajectory in the
time interval $[B,B+T]$.  This parameter is necessary since if we
recorded every single state that the trajectory passed through, we
could quickly run out of memory.

Choosing sensible values for these parameters is not easy, since it is
problem dependent. However, our aim is not to run such a long
trajectory that simply using a histogram of the values would lead to
an accurate approximation of the invariant density, as this would
completely negate the need for the rest of the algorithm.
Likewise, if the trajectory is too short, we will have a very poor
approximation of the areas of the domain which contain the majority of
the probability density. Our aim, therefore, is to get as good an
approximation of this region as possible with the shortest possible
trajectories.

As with many Monte Carlo estimates, due to the law of large numbers,
the error decays at a rate proportional to $1/\sqrt{T}$. It is
reasonable to expect that the approximation of the region which
contains the majority of the invariant density
should decay at a similar rate. Since convergence diagnostics relating
to Markov chains are an open area with no hard and fast
solution\cite{Cowles:1996:MCM}, we suggest that a sample trajectory be
created a priori to ascertain the relaxation time of the system. The
parameter $B$ can simply be set to be a proportion of the length $T$,
for instance $B = T/10$. The parameter $Q$ should be picked according
to how much memory one would like to set aside to store the sampled
points. Since we are calculating minimum distances to any point in
this set in step {\bf [6]} of the algorithm, the more points there
are, the longer the algorithm will take to run.

  \subsection{Domain and mesh generation parameters}
  The domain and mesh generation parameters are given by $\beta_1 > 0$,
  $\beta_2 > 0$ and $H\in \mathbb{N}$. The parameter $\beta_2$ defines the
  size of the domain as seen in \eqref{eq:domain}. Unless the trajectories
  from the SSA are very short, the estimation of the domain should be
  relatively trivial, and a value of order 1 should suffice
  in most instances. The value of $\beta_1$, which indicates how much
  error we estimate there to be in our estimation of the region which
  contains the majority of the invariant density,
  should be directly proportional to $1/\sqrt{T}$. Too large a value
  will lead to unnecessary refinement in some areas, while too small a
  value will lead to not enough refinement and more error.
  The parameter $H$ indicates the maximum number of splits that the
  initial cubes of the mesh may undergo. In practice this can be found
  at runtime, by capping the total number of elements we allow in the
  mesh due to memory and CPU constraints.

  \subsection{A numerical test of accuracy}
  \label{sectionaccuracy}
  Although the method has seven parameters which have to be specified,
  it is relatively easy to check that the computed results are
  numerically accurate.  Given the parameter set $\{ S, T, Q, B,
  \beta_1, \beta_2, H \}$, we can compare the computed results with
  the results obtained with the parameter set $\{ S, 2T,
  2Q, 2B,
  2\beta_1, \beta_2, H+1 \}$. If the result does not change
  significantly, then we can conclude that our parameters were chosen
  sufficiently large. Here, the parameters $T$, $Q$, $B$, and
  $\beta_1$, are multiplied by $2$ because increasing any of these
  parameters will increase the accuracy of the solution. In a similar
  way $H$ can be increased by one, which means that the finest level
  of refinement in the mesh becomes smaller by a factor of two in each
  coordinate. However, we do not propose to modify $S$ and $\beta_2$
  for the following reasons. The number of starting points $S$ can be
  determined by the number of stable equilibria of (\ref{eq:ode}). In
  this case, the parameter $S$ does not have to be changed during this
  a posteriori test of accuracy. Furthermore, $\beta_2$ is used only
  to choose the size of the domain, and multiplying the size of the
  domain by $2$ would lead to a smaller proportion of the domain being
  covered by $\Gamma$ (given by \eqref{definitionofgamma}), and as
  such a more accurate solution would not be guaranteed. If one wishes
  to test whether $\beta_2$ is large enough, the parameter $H$ may
  have to be increased as a function of $\beta_2$ in order to maintain
  the same level of refinement in $\Gamma$.

\section{Implementation of the saFEM}
\label{sec:Imp}

As soon as the computational mesh is determined, the approximation
$p_h \in V_h$ of the invariant distribution $p$ is computed by the
standard FEM approach \cite{Ciarlet:1978:FEM,Solin:2008:ALH}.
In order to compute the entries of the stiffness matrix (\ref{eq:stiff}),
it is necessary to use suitable quadrature
routines. In numerical examples below, we consider chemical reactions
of at most second order. Therefore, the diffusion  $D(\vecx)$
and drift $\vecv(\vecx)$ coefficients are polynomials of degree at most two.
Since the finite element space consists of piecewise linear functions,
it suffices to use a tetrahedral quadrature rule of order three which
integrates cubic polynomials exactly. The optimal (Gauss) quadrature rule
on tetrahedra of order three has 5 points \cite{Solin:2008:ALH}.

Since the dimension of the resulting algebraic system
(\ref{eq:linsys}) can be very large, even for relatively coarse
approximations, parallelisation of the assembly process is of
paramount importance. Several packages exist for constructing matrices
in parallel. The Portable Extensible Toolkit for Scientific
Computation (PETSc) is a suite of data structures and routines for the
scalable (parallel) solution of scientific applications modelled by
PDEs\cite{Balay:1997:MST,Balay:2011:PWP}. The matrix is split into
sections which are controlled by each processor. If the degrees of
freedom are ordered in a sensible way, then Message Passing Interface
(MPI) traffic between the processors can be kept to a minimum, leading
to excellent scaling of processors vs. computation
time\cite{Balay:1997:MST}.

Once the stiffness matrix (\ref{eq:stiff}) is assembled, we have to
find a nontrivial solution of the algebraic system (\ref{eq:linsys}).
Practically, we look for the eigenvector corresponding to the zero
eigenvalue of the matrix $A$.  To find this eigenvector in parallel,
we use a sister package of PETSc, which is called the Scalable Library
for Eigenvalue Problem Computations (SLEPc)
\cite{Hernandez:2005:SSF}. In particular, we used the MUltifrontal
Massively Parallel sparse direct Solver (MUMPS)
\cite{Amestoy:2001:FAM,Amestoy:2006:HSP} package for the
preconditioning, and a power method to solve the resulting eigenvalue
problem \cite{Kubla:1961:SCE}.

Once we have calculated the eigenvector, it is then necessary to
reconstruct the approximated function. For this, we use
the information regarding the hanging nodes so that we reconstruct
all of the vertices on the mesh. The function in question is a
probability density and therefore we normalise the obtained finite
element solution $p_h$ such that it satisfies the condition
(\ref{eq:strong1}).

\section{Numerical Results}
\label{sec:res}

In this section, we first use a simple example chemical system from
\cite{Cotter:2011:CAM} and implement the saFEM on a range of different
mesh sizes, and with different algorithmic parameters. Then we present
the results of saFEM for the Oregonator \cite{Field:1974:OCS}.

\subsection{Convergence of the numerical method}

We will study the system of three chemical
species $X_1$, $X_2$ and $X_3$ which are subject to the
following system of five chemical reactions \cite{Cotter:2011:CAM}:
\begin{equation}\label{eq:lin}
\RARR{\RARR{\RARR{\emptyset}{\LRARR{X_1}{X_2}{k_{2}}{k_{3}}}{k_1}}{X_3}{k_4}}{\emptyset}{k_5}.
\end{equation}
Then the propensity functions are defined by
\begin{eqnarray*}
&&\alpha_1(t) = k_1 V, \quad
\alpha_2(t) = k_2 X_1(t), \quad
\alpha_3(t) = k_3 X_2(t),
\\
&&\hspace{12mm}\alpha_4(t) = k_4 X_2(t), \quad
\alpha_5(t) = k_5 X_3(t),
\end{eqnarray*}
where $V$ is the volume of the reactor
\cite{Cotter:2011:CAM,Erban:2007:PGS}.
We consider this system with the following set
of non-dimensionalised parameters:
\begin{equation}\label{eq:linparam}
k_1V = 100, \quad
k_2 =  k_3 = 5, \quad
k_4 = k_5 = 1.
\end{equation}

As discussed in Section \ref{sec:param}, there are seven different
algorithmic parameters, whose values determine the accuracy and
efficiency of the methods. In this section we will analyse the effects
and convergence of the method due to altering the three most important
of these parameters.

\subsection{Convergence of approximation}
First we consider how the error in the approximation of the invariant
distribution $p$ decays as we refine the mesh, each time using the
same stochastic simulation, with $S=1$, $T = 10^5$, $Q=0.1$, $B=10^3$,
$\beta_1 = 0.01$ and $\beta_2 = 0.55$. We choose $S=1$ since the
corresponding ODE approximation given by \eqref{eq:ode} gives us:
\begin{eqnarray*}
\frac{\mbox{d}x_1}{\dt} &=& 100 - 5x_1 + 5x_2, \\
\frac{\mbox{d}x_2}{\dt} &=& 5x_1 - 6x_2 - 0.5,\\
\frac{\mbox{d}x_3}{\dt} &=& x_2 - x_3,
\end{eqnarray*}
which has a single steady state at $(119.5, 99.5, 99.5)$. Therefore we
pick the single starting point for the SSA trajectory to be at $(120,
100, 100)$. Note that the constant term $-0.5$ in the second equation
comes from the derivative of the diffusion terms as given by
\eqref{eq:CFPE2}, while these derivatives sum to zero for the first
and third equations.

\begin{figure}[t]
\picturesAB{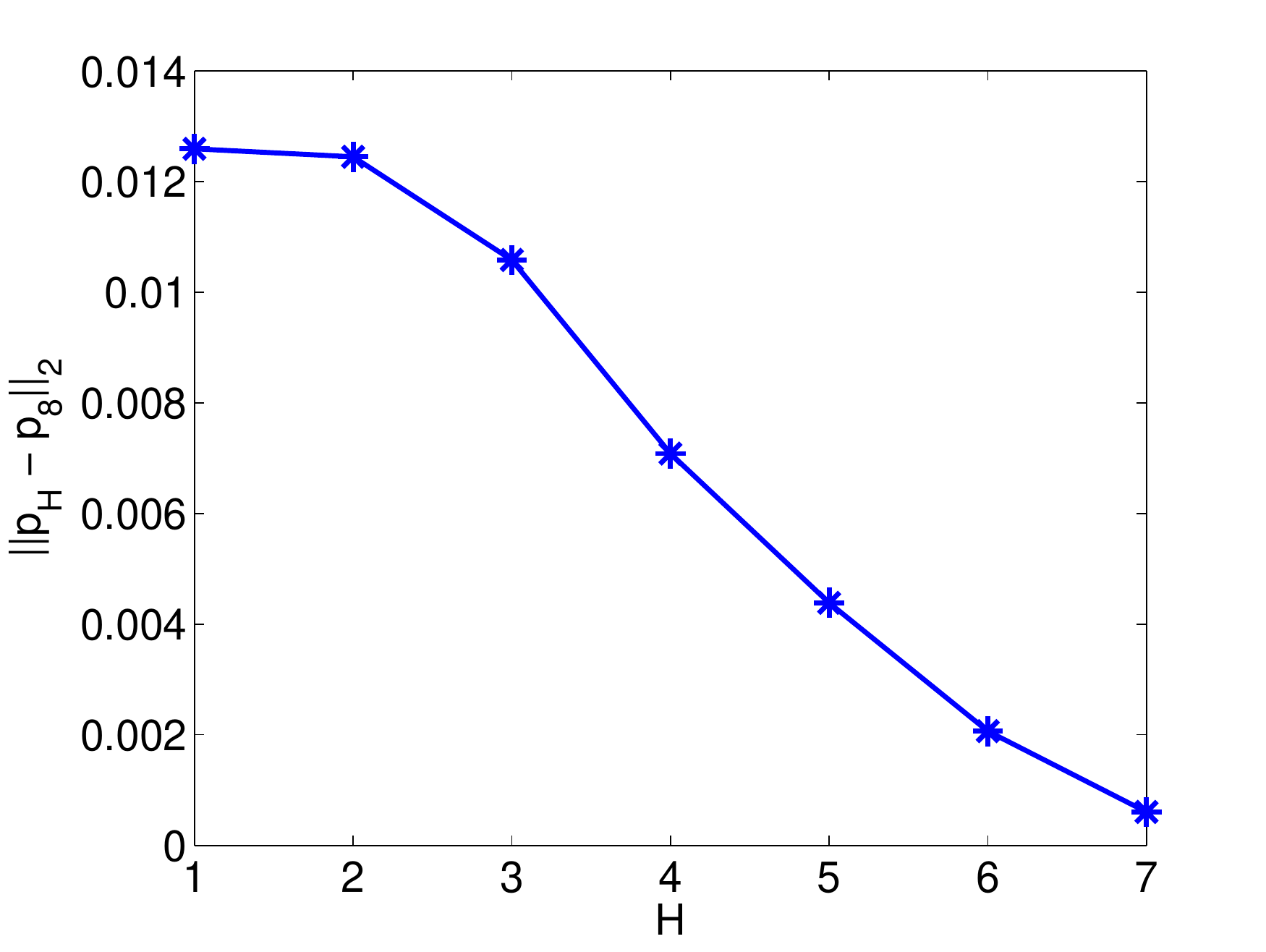}{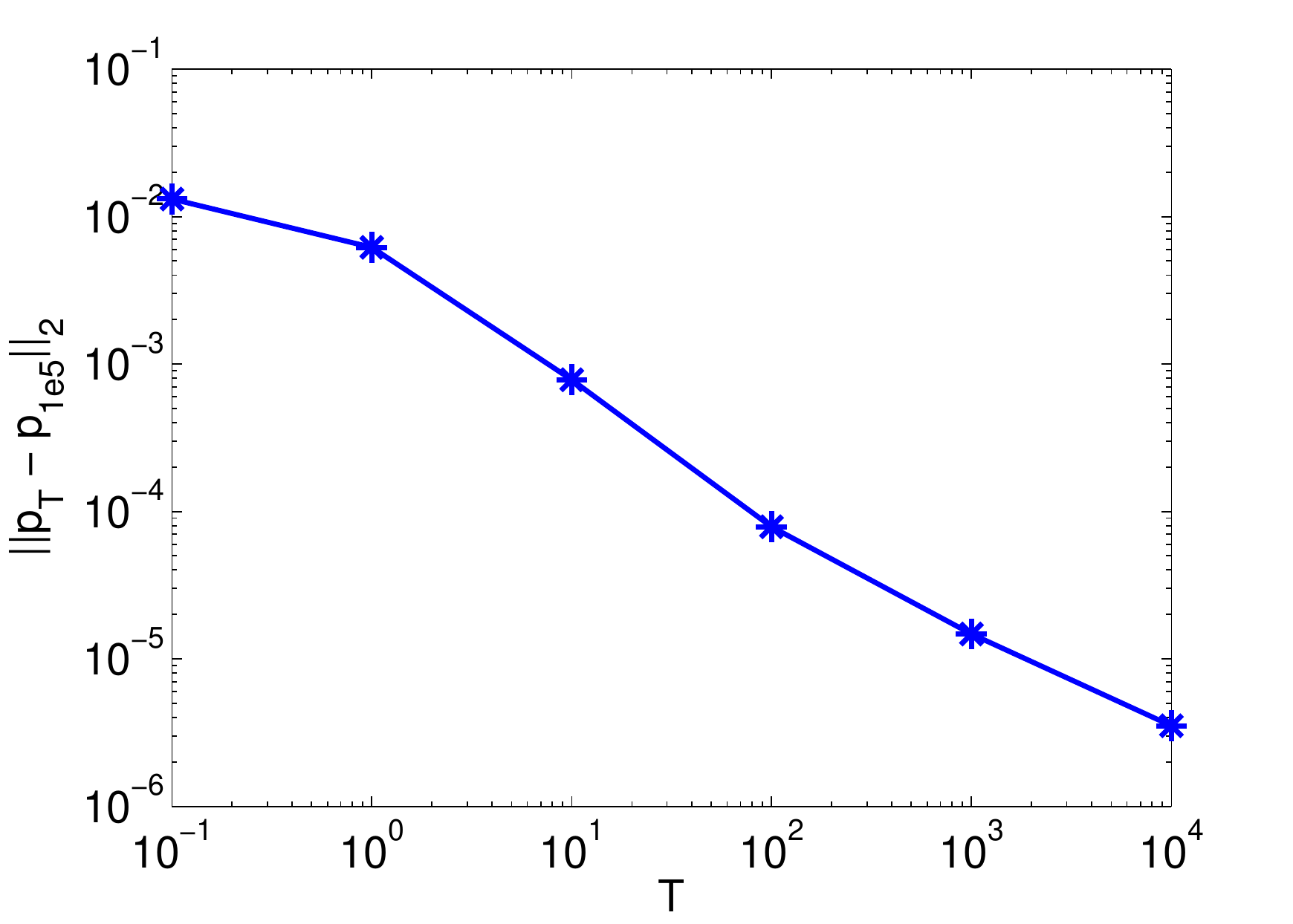}
{1.8in}{3mm}{-5mm}{5mm}
\caption{{\rm (a)} {\it Convergence of the saFEM with
$S=1$, $T = 10^5$, $Q=0.1$, $B=10^3$, $\beta_1 = 0.01$ and
$\beta_2 = 0.55$ over a range of mesh refinements $H$,
for the system $\eqref{eq:lin}$. Errors are given
by $(\ref{errorformula})$.}
{\rm (b)} {\it Convergence of the saFEM with $S=1$,
$Q=10$, $B=10^3$, $\beta_1 = 0.01$, $\beta_2 = 0.4$,
and $H=6$ over a range of SSA lengths $T$, for the
system \eqref{eq:lin}. Error is given by $(\ref{errorformula2})$.}
\label{fig:convab}}
\end{figure}

Figure \ref{fig:convab}(a) shows the convergence of the method as the
number of splits $H$ is increased from 1 to 7, where the error is
calculated by comparison with an approximation using a mesh which has
undergone up to 8 splits:
\begin{equation}
  \mbox{Error}(H) = \left (\int_\Omega |p_{8}({\bf x}) - p_{H}({\bf x})|^2\dx\right )^{1/2}.
\label{errorformula}
\end{equation}
When using a relatively coarse mesh, it is possible for the value of
the approximation to be negative in regions. Since the function
represents a probability density, which is strictly non-negative, this
does not make sense. These negative areas can also complicate the
normalisation of the function. We solve this problem by simply cutting
off the negative values, i.e. we multiply the approximated function by
the indicator function over the set where the function is
non-negative, before normalising by (\ref{eq:strong1}).

\subsection{Length of Stochastic Simulation}
Next we show convergence of the method as $T$ is increased, with
$S=1,$ $Q=10$, $B=10^3$, $\beta_1 = 0.01$, $\beta_2 = 0.4$ and $H=8$.
As in the previous section, the SSA trajectory is given initial
condition $(X_1, X_2, X_3) = (120, 100, 100)$.

As the system is ergodic, we would expect that as the length of
stochastic simulation
increases and the number of samples increases, that our samples
becomes increasingly representative of the invariant density in
question. This gives us more information about where we should be
refining our mesh, which in turn should give us a more accurate
approximation.

Figure \ref{fig:convab}(b) shows the convergence of approximation as
the length of stochastic simulation is increased. The error as plotted
in Figure \ref{fig:convab} (b) is the $L^2$
difference between the approximations with varying $T$ and the
approximation calculated with $T=10^5$, namely:
\begin{equation}
\mbox{Error}(T) = \left (\int_\Omega |p_{10^5}({\bf x}) - p_{T}({\bf x})|^2\dx\right )^{1/2}.
\label{errorformula2}
\end{equation}
Note that in order to more easily compare the distributions, the
domain selection for all of the data points was made using the longest
stochastic simulation with $T=10^5$.  This shows that as the
simulation length is increased we see a convergence of the chosen mesh
to one which well represents the region containing the invariant
density.

\subsection{Parameter $\beta_1$}
\begin{figure}[hb]
\centering
\includegraphics[width=0.45\textwidth]{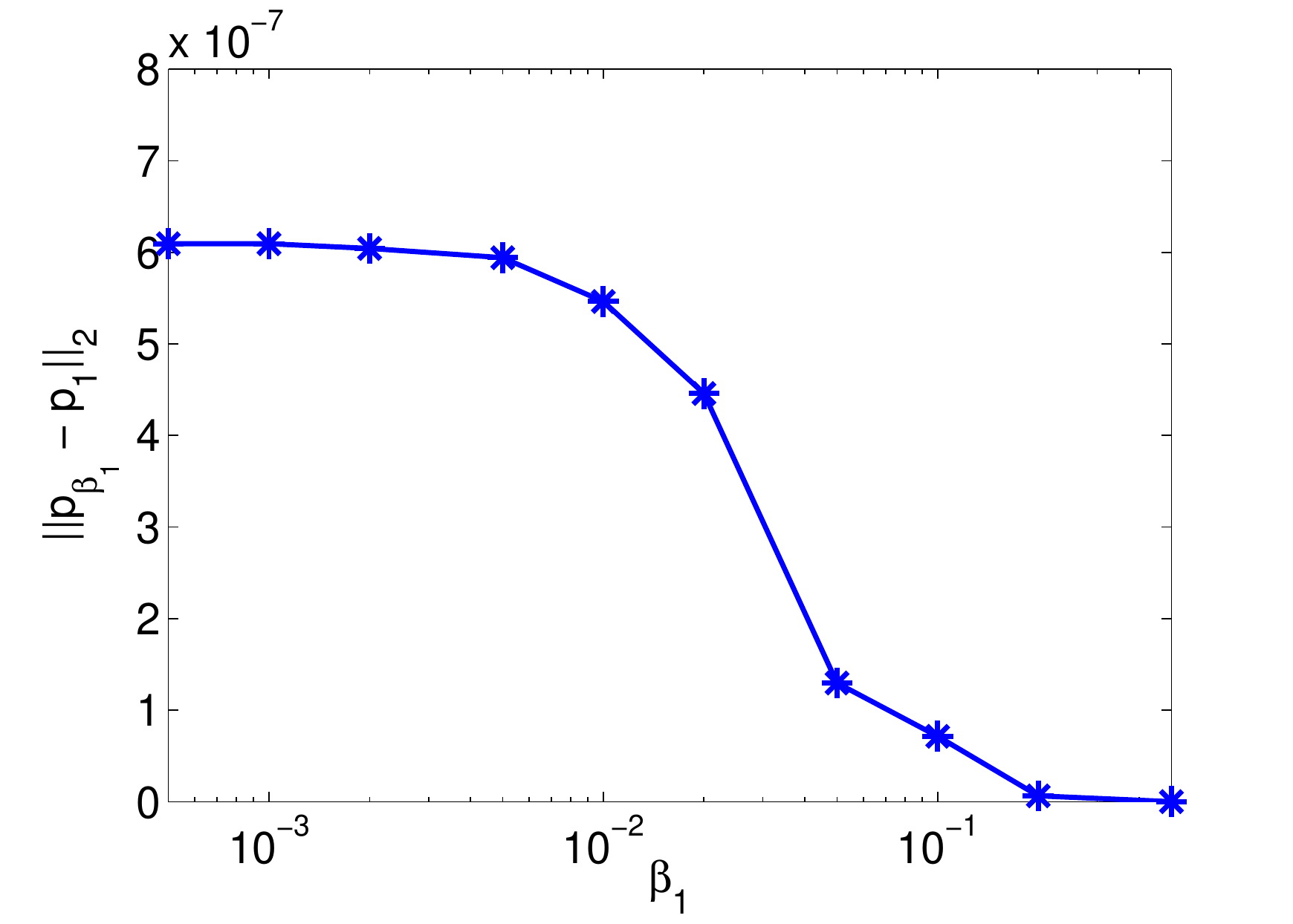}
\caption{\it Convergence of the saFEM with $\beta_2 = 0.4$, $Q=0.1$,
  $B=10^3$, $T=10^5$ and $S=6$ over a range of radii of uncertainty
  $\beta_1$, for the system \eqref{eq:lin}. The error is defined by
  \eqref{errorformula3}. \label{fig:r}}
\end{figure}

In certain situations the stochastic simulations may be expensive and
we may only be able to take a few samples from the stochastic
trajectory. In this case we still have a large amount of uncertainty
about the region which contains the vast majority of the invariant
density. Therefore we can only hope to approximate this by increasing
the size of the ellipsoid around each sampled point that we include in
the region $\Gamma$ given by \eqref{definitionofgamma}.  Figure
\ref{fig:r} shows the convergence of the method as $\beta_1$ is
increased, with $S=1$, $T=10^5$, $Q=0.1$, $B=10^3$, $\beta_2 = 0.4$
and $H=6$.

The error as plotted in Figure \ref{fig:r} is given by:
\begin{equation}
\mbox{Error}(\beta_1) = \left (\int_\Omega |p_{1}({\bf x}) - p_{\beta_1}({\bf x})|^2\dx\right )^{1/2}.
\label{errorformula3}
\end{equation}

\subsection{A numerical test of accuracy}
In Section \ref{sectionaccuracy} a brief test for sufficient
convergence was suggested, where two solutions were calculated, one
with parameters 
\[\{ S, T, Q, B, \beta_1, \beta_2, H \},\]
and the second with parameters \begin{equation}\label{P2} \{ S,
  2T, 2Q, 2B,
2\beta_1, \beta_2, H+1 \}.\end{equation} If the two
solutions are close up to some tolerance, then we can be fairly sure that
the method has converged to a reasonable degree.

The solution for the system \eqref{eq:lin} calculated with the
parameter set \[\{ S=1,\; T=10^5,\; Q=0.05,\; B=500,\; \beta_1=10^{-3},\; \beta_2
=0.45,\; H=7 \},\] when compared with the solution using the parameters
\[\{S=1,\; T=2\times 10^5,\; Q=0.1,\; B=10^3,\; \beta_1=2\times
10^{-3},\; H=8 \},\]  gave $L^2$ error of $1.99\times 10^{-4}$, which
is significantly less than the $L^2$ norm of the solution ($4.51
\times 10^{-3}$).  This gives a relative $L^2$ error of $4.41\times
10^{-2}$. For ease of comparison the domain $\Omega$ chosen
using parameters \eqref{P2} was used for both meshes.

\subsection{Oregonator}
Our final example is the Oregonator\cite{Field:1974:OCS} which is
a three species chemical system motivated by the Belousov-Zhabotinsky
reaction. It exhibits oscillatory behaviour and is traditionally
given by the following system of ODEs:
\begin{eqnarray}
\frac{\mbox{d}x_1}{\dt} &=& k_1x_2 - k_2x_1x_2 + k_3x_1 - 2k_4x_1^2,
\nonumber\\
\frac{\mbox{d}x_2}{\dt} &=& -k_1x_2 - k_2x_1x_2 + \frac{1}{2}k_cx_3,
\label{eq:OregonODE}\\
\frac{\mbox{d}x_3}{\dt} &=& 2k_3x_1 - k_cx_3.\nonumber
\end{eqnarray}
We now construct a set of reactions whose behaviour is given by
\eqref{eq:OregonODE} in the thermodynamic limit: 
\begin{eqnarray}\label{eq:Oregon}
\RARR{X_2}{X_1}{k_1} \qquad & \RARR{X_1+X_2}{\emptyset}{k_2} &
\qquad \RARR{X_1}{2X_1+2X_3}{k_3} \\
\RARR{2X_1}{\emptyset}{k_4} &
\qquad \RARRlong{X_3}{\emptyset}{k_5 = k_c/2}
& \qquad \RARRlong{X_3}{X_2}{k_6 = k_c/2}. \nonumber
\end{eqnarray}
Note that the reaction with parameter $k_c$ has been split into two
reactions $R_5$ and $R_6$ in order to get the factor of half in the
equation for the rate of change of $x_2$, with $k_5 = k_6 = k_c/2$.
We consider this system with the following set of dimensionless
parameters:
\begin{equation}
\label{eq:Oregon:params}
k_1 = 0.3, \quad k_2 = 4000, \quad  k_3 = 5,\quad
k_4 = 1200, \quad k_c = 0.02.
\end{equation}
In this parameter regime the stochastic description of these reactions
exhibits oscillatory behaviour.
\begin{figure}[t]
\picturesAB{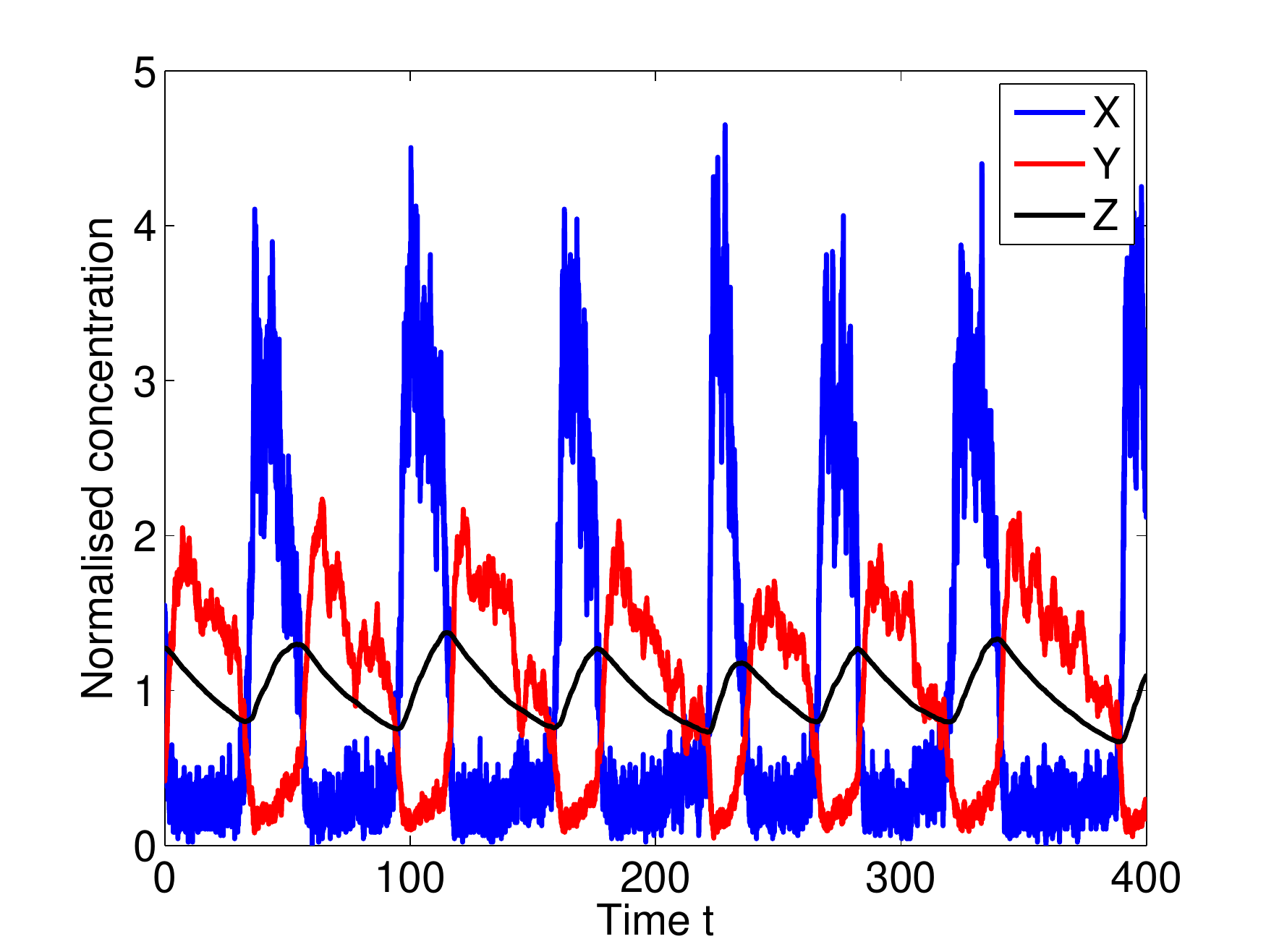}{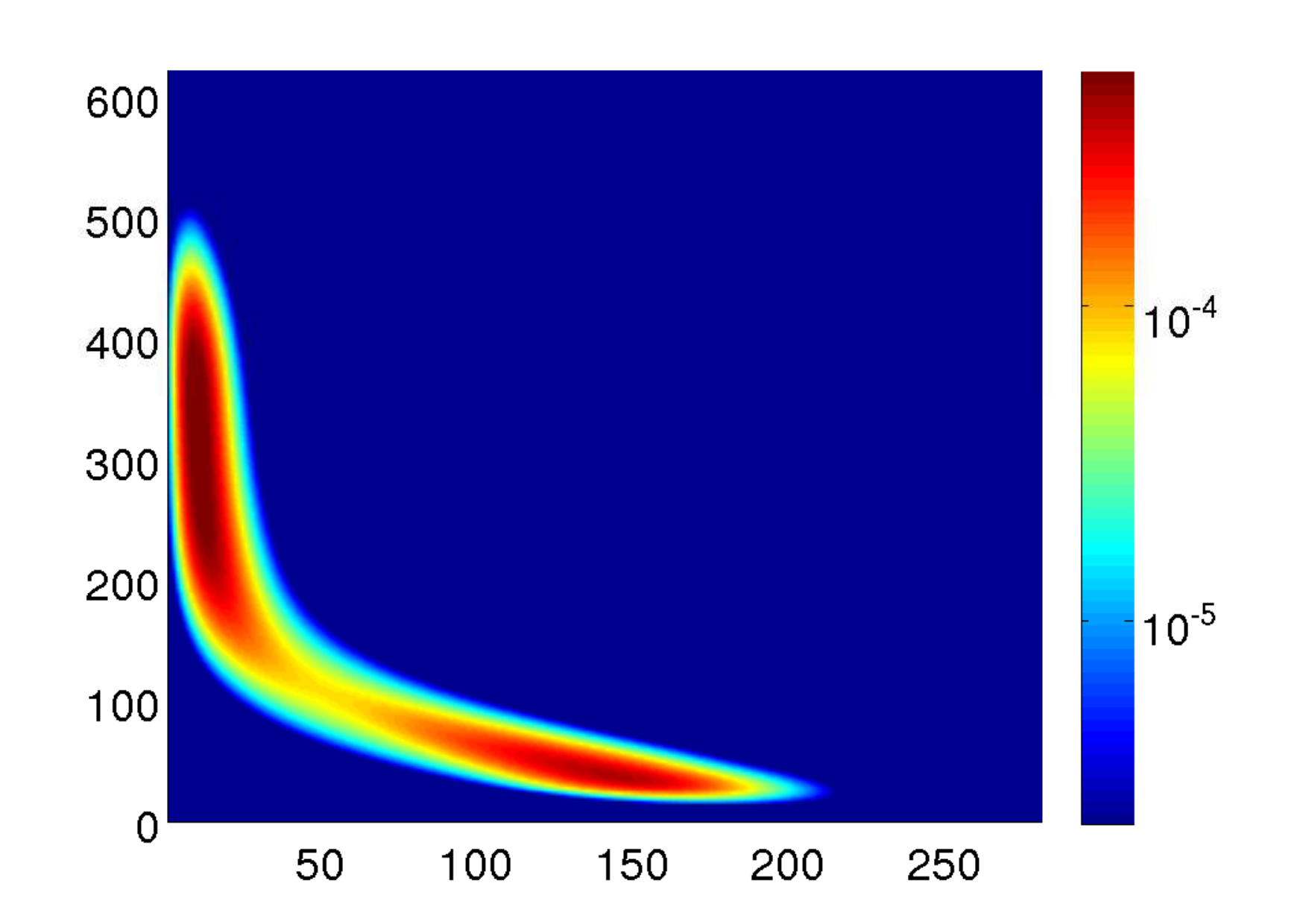}
{1.8in}{4mm}{-5mm}{8mm}
\vskip 3mm
\picturesCD{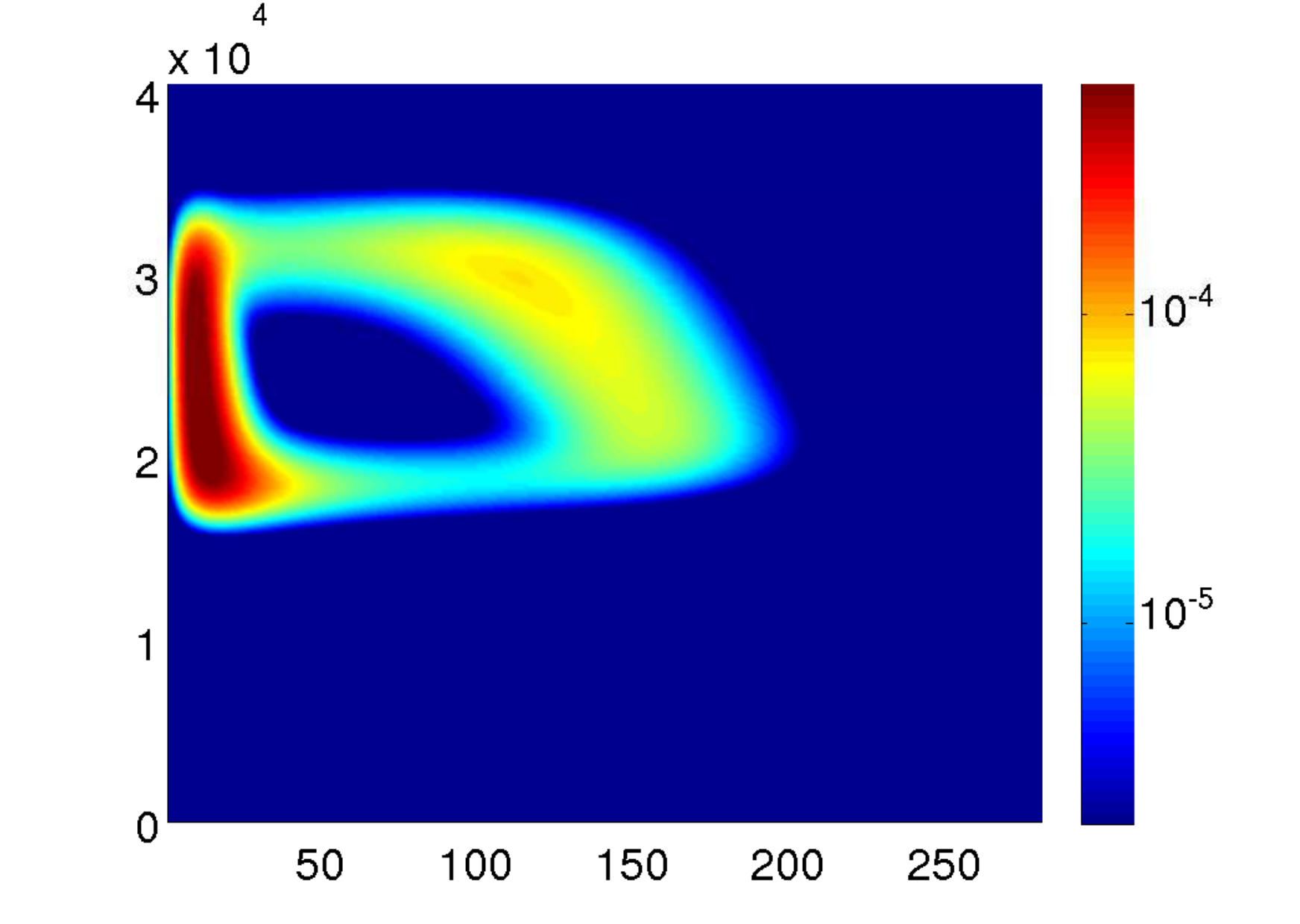}{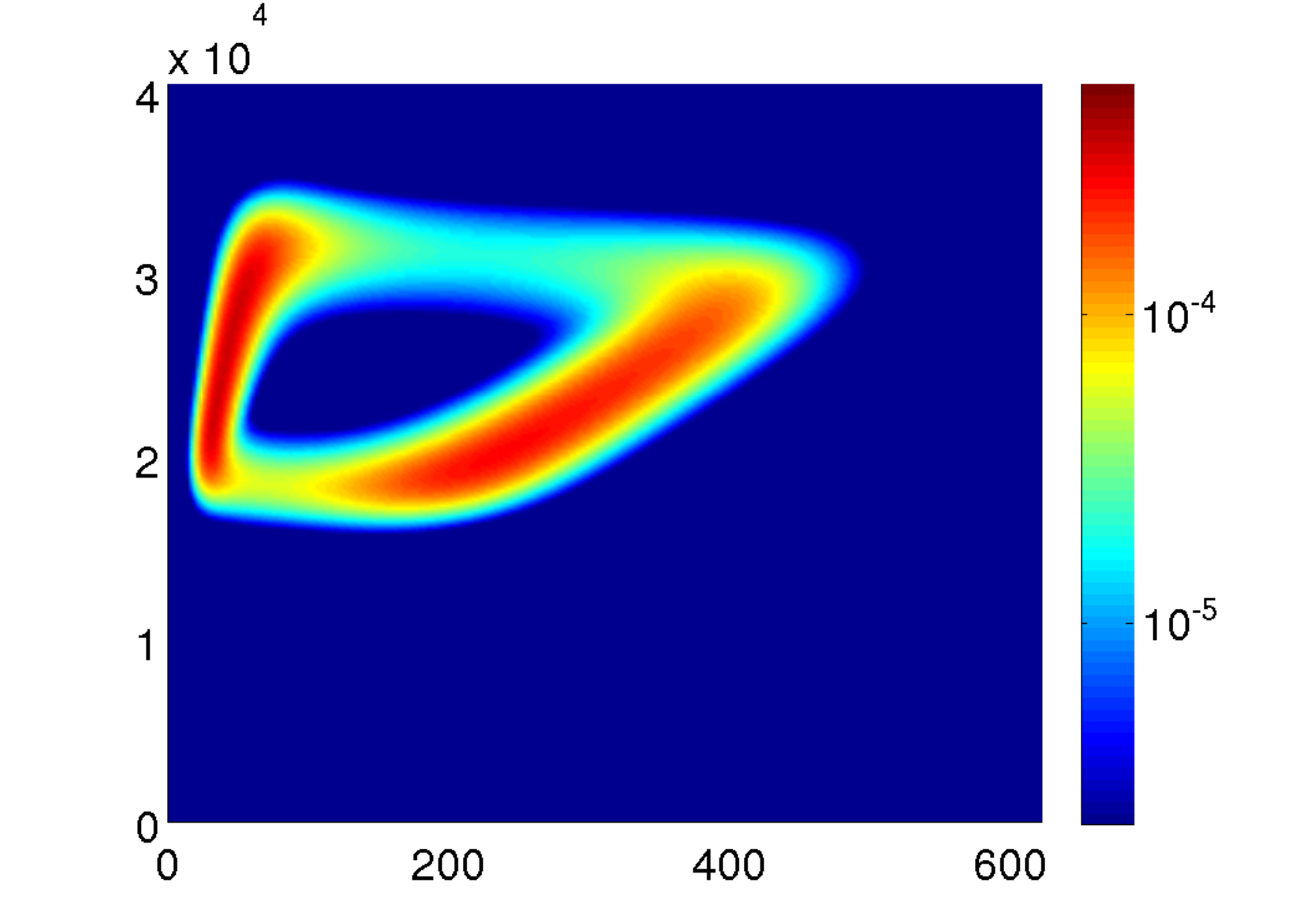}
{1.8in}{5mm}{-5mm}{5mm}
\vskip 2mm
\centering
\begin{picture}(0,0)
\put(18,260){$x_2$}
\put(103,167){$x_1$}
\put(20,110){$x_3$}
\put(103,15){$x_2$}
\put(-184,110){$x_3$}
\put(-100,15){$x_1$}
\end{picture}
\caption{{\rm (a)} {\it Trajectories of the Oregonator \eqref{eq:Oregon} with
  parameters given by \eqref{eq:Oregon:params}, simulated using the
  Gillespie SSA. Plots show normalised trajectories, with actual
  molecule numbers divided by a time average ($47.7$, $211.5$ and
  $2.4\times 10^4$ resp.). {\rm (b)--(d)} Approximation of the log of the
  marginal invariant distribution (conditioned on non-extinction of
  species $X_1$) of \eqref{eq:Oregon} with system parameters
  \eqref{eq:Oregon:params} in the {\rm (b)} $x_1$-$x_2$ plane,
  {\rm (c)} $x_1$-$x_3$ plane, {\rm (d)}
  $x_2$-$x_3$ plane, by the saFEM with algorithmic parameters given by
  \eqref{eq:params}}.\label{big_fig}}
\end{figure}
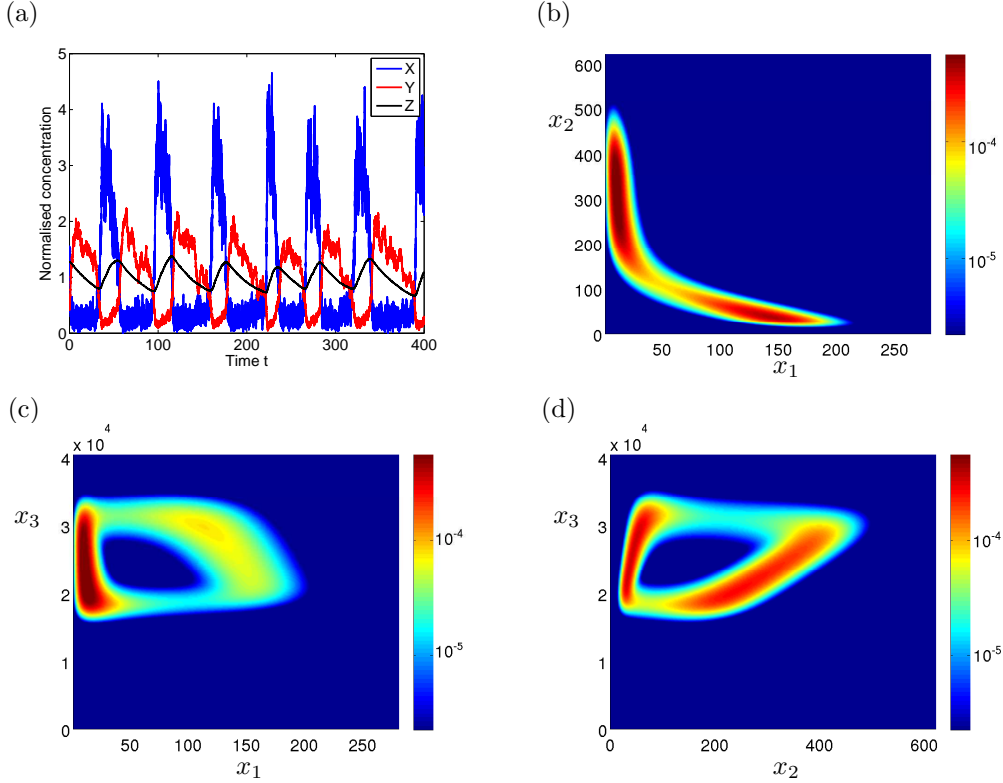
Figure \ref{big_fig}(a) shows normalised trajectories of the
Oregonator \eqref{eq:Oregon} with parameters given by
\eqref{eq:Oregon:params}, simulated using the Gillespie SSA. This
figure demonstrates the oscillatory behaviour of the Oregonator in
this parameter regime.

One thing of note should be mentioned at this point, regarding the
ergodicity of this system. The zero state, at the origin, is an
absorbing state for this system, and so trivially the invariant
distribution for this system is a Dirac measure on this
state. However, we are interested in the behaviour of this system
conditioned on non-extinction of the species. Therefore we ensure that
our domain of solution does not include this state, and thus the
transient states involved in the oscillation behaviour now form the
regions with non-zero invariant density.

In the figures that follow, these algorithmic parameters were used to
approximate the steady state distribution (conditioned on
non-extinction of species $X_1$):
\begin{equation}\label{eq:params}
  S=1, \qquad T=10^6, \qquad Q = 10^{-2}, \qquad B=10^3, \qquad  H = 8
\end{equation}
Since this system has a single limit cycle and is ergodic, step [{\bf 1}] of
the algorithm (in Table \ref{Meshgen Table}) can be replaced by
running the Gillespie SSA for a period of algorithm time until we are
reasonably sure that the chain has entered probabilistic equilibrium,
and using the endpoint of this simulation as the starting point for
the SSA in step [{\bf 2}]. Since the system has a single stable orbit, we
pick $S = 1$. Since a priori we do not know a specific point on the
most likely orbit, the initial condition $(x_1, x_2, x_3) =
(100,100,100)$ was chosen.

Since we wish to condition on non-extinction of all of the species, we
must ensure that $X_1$ never drops below 1. We can do this by slightly
altering the domain $\Omega$ of solution of the FPE. We do this by
replacing the first line of \eqref{eq:domain} by
$\Omega =  {\mathcal A}^*_1
\times {\mathcal A}_2
\times {\mathcal A}_3,$ where
\begin{equation*}
{\mathcal A}^*_1
=
\left( \max \{1, x_1^{\rm min}- \beta_2 \, x_1^{\rm range}\},
x_1^{\rm max} + \beta_2 \, x_1^{\rm range}
\right).
\end{equation*}
The longer the simulation in the step {\bf [2]} of the saFEM, the more
sure we can be of the regions which have non-negligible probability density,
and therefore the larger we can make $\beta_2$ and $\beta_1$. We
choose $\beta_1 = 0.01$ and $\beta_2=0.1$.

Using these parameters, a mesh was created by the saFEM with $8.65\times 10^5$
vertices, out of a possible $1.70\times 10^7$ with up to 8 splits of each
cuboid, with $1.17 \times 10^5$ hanging nodes. Therefore there were a total
of $7.48\times 10^5$ degrees of freedom, a mere $4.4\%$ of the degrees of
freedom that would have been required to fill the domain with cuboids of the
finest refinement level. The mesh contained $4.79\times 10^6$ individual
tetrahedral elements.

\begin{figure}[t]
\centering
\includegraphics[width=0.8\textwidth]{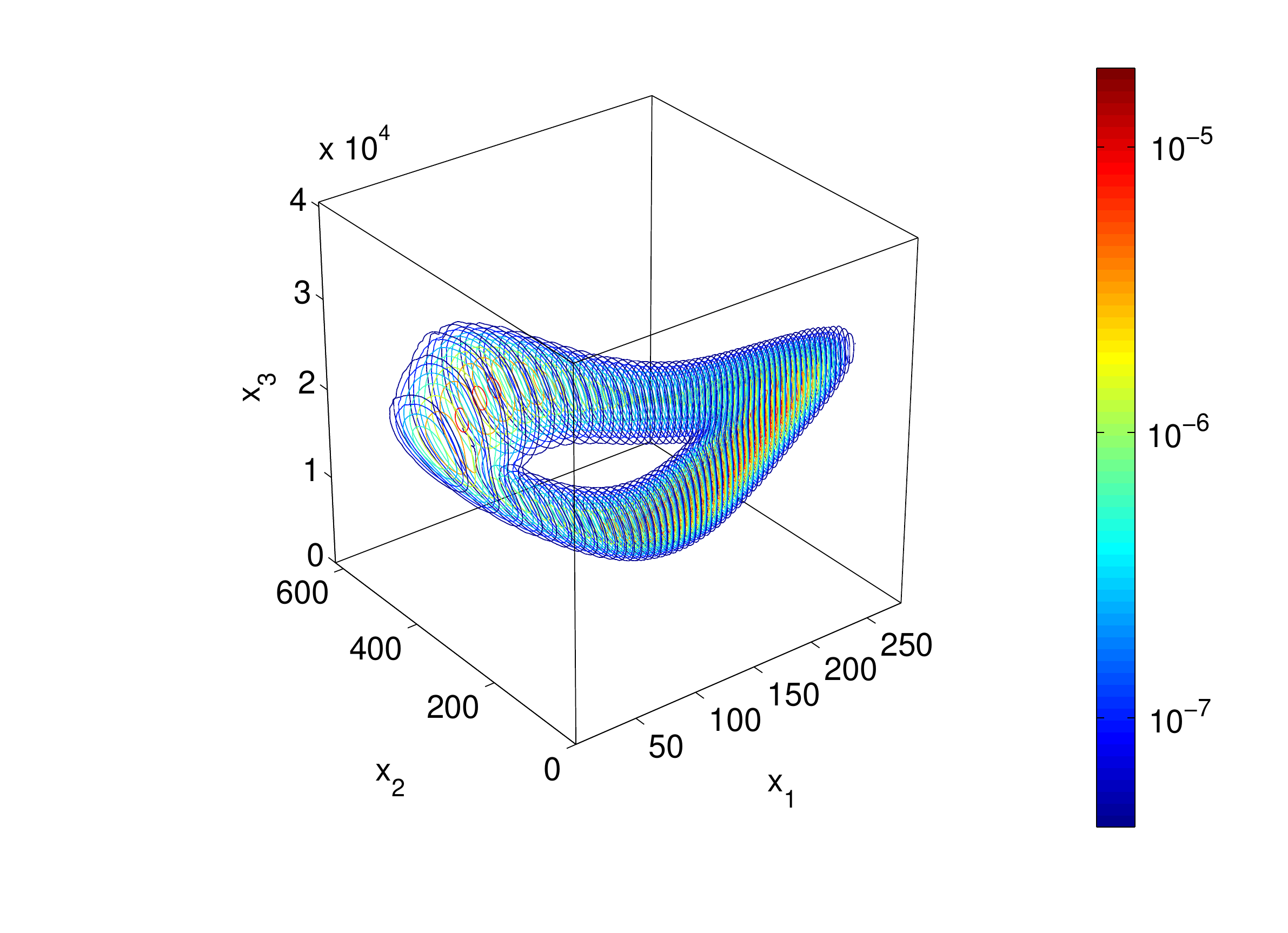}
\caption{\it Approximation of the invariant distribution (conditioned
  on non-extinction) of \eqref{eq:Oregon} with system parameters
  \eqref{eq:Oregon:params}. Contours plots shown in several
  slices.\label{contours}}
\end{figure}
Figure \ref{contours} shows several slices of the approximated
invariant probability density conditioned on non-extinction of species
$X_1$, represented by contour plots.
Since visualisation of a three-dimensional function on a
two-dimensional page is problematic, we present also in Figure
\ref{big_fig}(b)--(d) the marginal densities of the approximation of
the invariant density of the Oregonator system with parameters
\eqref{eq:Oregon:params} in 3 different planes.

To verify that the method is accurately representing the invariant
distribution, we can test the approximation as suggested in Section
\ref{sectionaccuracy}. The approximation was compared with one with
the parameters given as follows:
\begin{equation*}
  S=1, \; T=5\times 10^5, \; Q = 5\times 10^{-3}, \;
  B=5\times 10^2, \;  H = 7, \; \beta_1 = 5\times 10^{-3}.
\end{equation*}
We omit $\beta_2$ here as the same domain as was chosen for the
original parameter set was used for ease of comparison of the two
densities. The $L^2$ difference between the two solutions was
$2.89\times 10^{-11}$, with the original solution having an $L^2$
norm of $1.74 \times 10^{-8}$, giving a relative $L^2$ error of
$1.66 \times 10^{-3}$, verifying that our solution is
well converged.
\section{Discussion}

In this paper we have presented a method for solution of the
steady-state Fokker-Planck equation for chemical systems. The method
is based on the use of stochastic trajectories to estimate the region
in which we must refine our finite element mesh the most. This method
is only valid for systems which do not have large regions of the
domain with non-negligible invariant density that can be accurately
approximated by a linear function. Our assumption that regions which
have high enough probability density require high levels of refinement
would fail in this case, and the mesh produced would be
inefficient.

The use of the Fokker-Planck equation as a means to model chemical
reactions of the mesoscale has been present in the literature for some
time\cite{Gillespie:1980:AME,Sjoberg:2009:FPA}. As computational
resources have improved, the solving of such equations in more than
one dimension has become tractable. It has been shown that using
appropriate methods (such as finite volume methods), the solution of
the Fokker-Planck equation is far more efficient than using an SSA for
the computation of both time-dependent and steady state solutions, and
with a high degree of accuracy\cite{Ferm:2006:CSF,Sjoberg:2009:FPA}.

The Fokker-Planck equation does suffer badly from the curse of
dimensionality, and many different methods exist which attempt to find
approximations of the solution of these equations efficiently in more
than one dimension. Some involve dimension reduction of the equations
themselves, through simplifying assumptions, and through coupling with
macroscopic reaction rate equations for some of the
species\cite{Lotstedt:2006:DRF}. Sparse grid methods have also been
used in an attempt to find the invariant density of the chemical
master equation\cite{Hegland:2008:SGH}. Sparse grids attempt to
overcome the curse of dimensionality by reducing greatly the number of
degrees of freedom considered. This CME solution can also be coupled
to Fokker-Planck equations for the more abundant species using a
hybrid method\cite{Hegland:2008:SGH}. Other methods have been used to
approximate the solution of the chemical master equation, including
adaptive wavelet compression\cite{Jahnke:2010:AWC} and finite state
projection\cite{Munsky:2006:FSP,Drawert:2010:DFS}.

There are several other ways that stochastic trajectories could be
calculated in the saFEM, other than the standard SSA as detailed in
Table \ref{SSAtable}. The $\tau$-leap method\cite{Gillespie:2001:AAS}
approximates the trajectory by modelling several reactions in each
iteration of the algorithm, up to a fixed time step. This can
accelerate the simulation of the trajectory, but could incur some
errors, including the possibility that the trajectory may become
negative in one or more of the chemical species. Likewise, one could
use numerical approximations of the chemical Langevin
equation\cite{Gillespie:2000:CLE}, a stochastic differential equation
(SDE) with corresponding Fokker-Planck equation given by
\eqref{eq:CFPE}, to create trajectories. As with the $\tau$-leap
method, the trajectories from this SDE can become negative, and so
with both of these methods proper treatment of boundary conditions is
key to accurately approximating the region of the positive quadrant
that contains the majority of the probability density.

More information could also be extracted from the stochastic
trajectories. A rudimentary approximation of the density could be made
using the samples taken from the simulations, and used as a
preconditioner for the eigensolver. As the eigensolver used is
iterative, even a rough guess of the form of the solution could help
to reduce computation times. It has been previously shown that
combining stochastic simulations and solutions of differential
equations can be used for preconditioning of
computations\cite{Qiao:2006:SDS}.

\medskip

{\bf Acknowledgements:} The research leading to these results has
received funding from the European Research Council under the {\it
  European Community}'s Seventh Framework Programme ({\it
  FP7/2007-2013})/ ERC {\it grant agreement} No. 239870. This
publication was based on work supported in part by Award No
KUK-C1-013-04, made by King Abdullah University of Science and
Technology (KAUST). The work of Simon Cotter was also partially
supported by a Junior Research Fellowship of St Cross College,
University of Oxford. Tom\'a\v s Vejchodsk\'y would like to
acknowledge the support by the Grant Agency of the Academy of Sciences
(project.\ IAA100760702) and the Academy of Sciences of the Czech
Republic (institutional research plan no.\ AV0Z0190503). Radek Erban
would also like to thank Somerville College, University of Oxford, for
a Fulford Junior Research Fellowship; Brasenose College, University of
Oxford, for a Nicholas Kurti Junior Fellowship; the Royal Society for
a University Research Fellowship; and the Leverhulme Trust for a
Philip Leverhulme Prize. This prize money was used to support a
research visit of Tom\'a\v s Vejchodsk\'y in Oxford.

\bibliography{bibrad}

\bibliographystyle{plain}

\end{document}